\documentclass[9pt, journal,comsoc]{IEEEtran}
\usepackage[T1]{fontenc}

\usepackage{cite}

\ifCLASSINFOpdf
\else
\fi
\usepackage{amsmath}
\usepackage{algpseudocode,algorithm2e}
\usepackage{hyperref}  
\usepackage{amsfonts}
\usepackage{mathtools}
\usepackage{subcaption}
\usepackage{amssymb}
\usepackage{bbm}
\usepackage{mathtools}
\usepackage{enumitem}

\usepackage{xcolor}
\usepackage{tikz}
\usepackage{tkz-graph}

\DeclareMathOperator{\spec}{spec}

\let\leq\leqslant

\let\emptyset\varnothing

\newcommand{\calA}{\ensuremath{\mathcal{A}}}
\newcommand{\calB}{\ensuremath{\mathcal{B}}}

\newcommand{\calD}{\ensuremath{\mathcal{D}}}

\newcommand{\calG}{\ensuremath{\mathcal{G}}}

\newcommand{\calI}{\ensuremath{\mathcal{I}}}

\newcommand{\calM}{\ensuremath{\mathcal{M}}}
\newcommand{\calN}{\ensuremath{\mathcal{N}}}

\newcommand{\calP}{\ensuremath{\mathcal{P}}}

\newcommand{\bbC}{\ensuremath{\mathbb{C}}}

\newcommand{\bbP}{\ensuremath{\mathbb{P}}}

\newcommand{\bbR}{\ensuremath{\mathbb{R}}}

\newcommand{\bmat}{\begin{matrix}}
\newcommand{\emat}{\end{matrix}}
\newcommand{\bbm}{\begin{bmatrix}}
\newcommand{\ebm}{\end{bmatrix}}
\newcommand{\bbma}{\begin{bmatrix*}[r]}
\newcommand{\ebma}{\end{bmatrix*}}
\newcommand{\bpm}{\begin{pmatrix}}
\newcommand{\epm}{\end{pmatrix}}
\newcommand{\bvm}{\begin{vmatrix}}
\newcommand{\evm}{\end{vmatrix}}
\newcommand{\bse}{\begin{subequations}}
\newcommand{\ese}{\end{subequations}}
\newcommand{\beq}{\begin{equation}}
\newcommand{\eeq}{\end{equation}}
\newcommand{\beqn}{\begin{equation*}}
\newcommand{\eeqn}{\end{equation*}}

\newcommand{\ben}{\renewcommand{\labelenumi}{\arabic{enumi}.}
\renewcommand{\theenumi}{\arabic{enumi}}\begin{enumerate}}

\newcommand{\een}{\end{enumerate}}

\newcommand{\beni}{\renewcommand{\labelenumi}{\roman{enumi}.}
\renewcommand{\theenumi}{\roman{enumi}}\begin{enumerate}}

\newcommand{\eeni}{\end{enumerate}}

\newcommand{\bena}{\renewcommand{\labelenumi}{\alph{enumi}.}
\renewcommand{\theenumi}{\alph{enumi}}\begin{enumerate}}

\newcommand{\eena}{\end{enumerate}}
\newcommand{\bit}{\begin{itemize}}
\newcommand{\eit}{\end{itemize}}

\newtheorem{thm}{Theorem}
\newtheorem{defn}[thm]{Definition}
\newtheorem{lem}[thm]{Lemma}

\newtheorem{cor}[thm]{Corollary}
\newtheorem{example}[thm]{Example}
\newtheorem{remark}[thm]{Remark}
\newtheorem{prb}[thm]{Problem}
\newtheorem{prp}[thm]{Proposition}

\newcommand{\bthm}{\begin{thm}}
\newcommand{\ethm}{\end{thm}}
\newcommand{\blem}{\begin{lem}}
\newcommand{\elem}{\end{lem}}
\newcommand{\bprop}{\begin{prp}}
\newcommand{\eprop}{\end{prp}}
\newcommand{\bex}{\begin{example}}
\newcommand{\eex}{\end{example}}
\newcommand{\bas}{\begin{assumption}}
\newcommand{\eas}{\end{assumption}}
\newcommand{\bre}{\begin{remark}}
\newcommand{\ere}{\end{remark}}
\newcommand{\bcor}{\begin{cor}}
\newcommand{\ecor}{\end{cor}}
\newcommand{\bdfn}{\begin{defn}}
\newcommand{\edfn}{\end{defn}}
\newcommand{\bcon}{\begin{conjecture}}
\newcommand{\econ}{\end{conjecture}}
\newcommand{\bali}{\begin{aligned}}
\newcommand{\eali}{\end{aligned}}
\newcommand{\eprb}{\end{prb}}
\newcommand{\bprb}{\begin{prb}}
\newcommand{\ecas}{\end{cases}}
\newcommand{\bcas}{\begin{cases}}

\newcommand{\nset}[1]{\ensuremath{\{1,2,\ldots,#1\}}}

\newcommand{\set}[2]{\ensuremath{\{#1\mid #2\}}}

\newcommand{\0}{\ensuremath{\mathbf 0}}

\newcommand{\BP}{\begin{IEEEproof}}
\newcommand{\EP}{\end{IEEEproof}}
\definecolor{myc1}{rgb}{1, 0.188, 0.188}
\definecolor{myc2}{rgb}{0.188,1, 0.188}
\definecolor{myc3}{rgb}{0.188, 0.188, 1}
\definecolor{myc4}{rgb}{0.458, 0.388, 0.678}
\definecolor{myc5}{rgb}{0.558, 0.188, 0.678}
\definecolor{myg1}{rgb}{0.858, 0.488, 0.678}
\definecolor{myg2}{rgb}{0.658, 0.588, 0.478}
\definecolor{myg3}{rgb}{0.458, 0.788, 0.478}
\definecolor{myg4}{rgb}{0.458, 0.588, 0.678}
\definecolor{myg5}{rgb}{0.458, 0.888, 0.778}

\tikzset{LabelStyle/.style = { rectangle, rounded corners, 
		minimum width = 10pt, fill = white,
		text = black }}

\usepackage{graphicx}  
\usepackage{xcolor}
\makeatletter
\def\mathcolor#1#{\@mathcolor{#1}}
\def\@mathcolor#1#2#3{%
	\protect\leavevmode
	\begingroup
	\color#1{#2}#3%
	\endgroup
}
\makeatother
\usepackage{hyperref}

\newcounter{todocounter}
\usepackage[colorinlistoftodos]{todonotes}

\begin{document}

\title{Strong Structural Controllability of Colored Structured Systems}
%
%
%

\author{Jiajia~Jia,
	Harry~L.~Trentelman,~\IEEEmembership{Fellow,~IEEE,},
	Nikolaos~Charalampidis
	and~M.~Kanat~Camlibel,~\IEEEmembership{Member,~IEEE}
	\thanks{The authors are with the Bernoulli Institute for Mathematics, Computer Science and Artificial Intelligence, University of Groningen,  9700 Groningen, The Netherlands (e-mail: {\footnotesize{\tt j.jia@rug.nl},
			{\tt h.l.trentelman@rug.nl}, {\tt charal.nik@gmail.com}, {\tt m.k.camlibel@rug.nl}).}}
}

%
%

\markboth{Journal of \LaTeX\ Class Files,~Vol.~14, No.~8, August~2015}%
{Shell \MakeLowercase{\textit{et al.}}: Bare Demo of IEEEtran.cls for IEEE Communications Society Journals}
%



\maketitle

\begin{abstract}
	    This paper deals with strong structural controllability of linear structured systems in which the system matrices are given by zero/nonzero/arbitrary pattern matrices.
	Instead of assuming that the nonzero and arbitrary entries of the system matrices can take their values completely independently, this paper allows equality constraints on these entries, in the sense that {\em a priori} given entries in the system matrices are restricted to take arbitrary but identical values. 
	To formalize this general class of structured systems, we introduce the concepts of colored pattern matrices and colored structured systems.
	The main contribution of this paper is that it generalizes both the classical results on strong structural controllability of structured systems as well as recent results on controllability of systems defined on colored graphs. 
	In this paper, we will establish both algebraic and graph-theoretic conditions for strong structural controllability of this more general class of structured systems.
\end{abstract}

\begin{IEEEkeywords}
	Strong structural controllability, Network analysis, Linear systems, Graph theory
\end{IEEEkeywords}

%
\IEEEpeerreviewmaketitle
	
	\section{Introduction}
	Controllability is a fundamental concept in systems and control.
	For linear time-invariant (LTI) systems of the form
	$$\dot{x} = Ax + Bu,$$
	controllability can be verified using the Kalman rank test or the Hautus test \cite{TSH2012}.
	Often, the exact values of the entries in the matrices $A$ and $B$ are not known, but only the underlying interconnection structure between the input and state variables is known exactly.
	In order to formalize this, Mayeda and Yamada \cite{MY1979} have introduced a framework in which, instead of a fixed pair of real matrices, only the so-called {\em zero/nonzero} structure of $A$ and $B$ is given.
	This means that each entry of these matrices is known to be either a {\em fixed zero} or an {\em arbitrary nonzero real number}.
	In addition, to guarantee the controllability of all possible LTI systems with such a given zero/nonzero structure, in \cite{MY1979} they introduced the concept of {\em strong structural controllability}.
	Since then, many contributions have been made on the topic of strong structural controllability. (See \cite{B1981, RSW1992, JFB2011, CM2013, MZC2014, TD2015} and the references therein.)
	
	Roughly speaking, two basic assumptions prevail in the aforementioned literature: (1) each entry of the system matrices is either a fixed zero or an arbitrary nonzero value, and (2) the nonzero entries take their values independently.
	Concerning the first of these assumptions, however, in many practical scenarios there might also be entries that can take arbitrary zero or nonzero values.
	Examples can be found in \cite{RM2009, MHM2019, PPKAI2019,JHHK2019} and the references therein. 
	In such scenarios, it is impossible to represent the system using a zero/nonzero structure.
	To deal with this, recently in \cite{JHHK2019} the notion of zero/nonzero structure has been extended to a more general \emph{zero/nonzero/arbitrary structure}, and thus a unifying framework for strong structural controllability was established. 
	Regarding the second of the above assumptions, we note that in physical systems it is often not the case that the nonzero entries in the system matrices can take their values independently.
	Indeed, some of the nonzero entries in the system matrices might have dependencies. (See \cite{MHM2017, TVDF2018, TDF2018, LM2019, JTBC20181, JTBC2018} and the references therein.) 
	In particular, in \cite{JTBC20181} and \cite{JTBC2018}, the situation was considered that prescribed nonzero entries in the system matrices are constrained to take identical (nonzero) values. These constraints can be caused by symmetry properties (\cite{MHM2017,TVDF2018}) or by physical constraints on the system \cite{LM2019}. 
	
	In this paper, we will explore the strong structural controllability of systems in which neither of the above two basic assumptions is satisfied. 
	More explicitly, the present paper will extend the approach taken in \cite{JTBC20181} and \cite{JTBC2018} using the newly introduced unifying framework from \cite{JHHK2019}. 
	That is, we will study strong structural controllability of systems in which the  zero/nonzero/arbitrary structure of the system matrices is given, and moreover, in which some of the entries in the system matrices are constrained to take identical values.
	Following the naming convention in \cite{JTBC20181}, \cite{JTBC2018} and \cite{JHHK2019}, we will call such kind of systems {\em colored structured systems}. 
	The main contributions of this paper are the following:
	\begin{enumerate}
		\item  We establish sufficient algebraic conditions for the strong structural controllability of a given colored structured system in terms of a full row rank test on two so-called {\em colored pattern matrices}.
		\item  We establish a  test for the full row rank property of a given colored pattern matrix in terms of a new color change rule and colorability of the graph associated with the pattern matrix. 
		In order to introduce this color change rule, we also establish a necessary and sufficient condition under which a given square colored pattern matrix is nonsingular.  
		\item  Based on the above results, we establish sufficient graph-theoretic conditions for strong structural controllability of colored structured systems.
	\end{enumerate}
	The outline of the paper is as follows.
	Section \ref{s:2} presents some preliminaries. 
	In Section \ref{s:3}, we formulate the problem treated in this paper in terms of colored structured systems. 
	In Section \ref{s:4}, we establish sufficient algebraic conditions for controllability of colored structured systems. We also provide a counterexample to show that these conditions are not necessary.
	Section  \ref{s:5} presents a necessary and sufficient graph-theoretic condition under which a given square colored pattern matrix is nonsingular.
	In Section \ref{s:6}, we introduce a new color change rule and the concept of colorability of the graph associated with a given colored pattern matrix. Based on these concepts, we establish a graph-theoretic condition under which a given colored pattern matrix has full row rank, and thus we obtain sufficient graph-theoretic conditions for strong structural controllability of colored structured systems.
	Finally, in Section \ref{S:conclusion}, we provide our conclusions. 
	
	\section{Preliminaries} \label{s:2}
	In this paper, we will use standard notation. 
	We denote by $\mathbb{C}$ and $\mathbb{R}$ the fields of complex and real numbers, respectively. 
	The spaces of $n$-dimensional real and complex vectors are denoted by $\mathbb{R}^{n}$ and 
	$\mathbb{C}^{n}$, respectively. 
	Likewise, the spaces of $n \times m$ real and complex matrices are denoted by $\mathbb{R}^{n \times m}$ and $\mathbb{C}^{n \times m}$.
	For a given $n \times m$ matrix $A$, the entry in the $i$th row and $j$th column is denoted by $A_{ij}$.  
	For a given square matrix $A$, we denote its determinant by $\det(A)$. 
	In addition, $I$ and $\0$ will denote identity and zero matrices of appropriate dimensions.
	For a given finite set $S$, its number of elements will be denoted by $|S|$.
	A finite collection $\{S_1, S_2, \ldots, S_k\}$ of subsets of $S$ is called a partition of $S$ if $S_i \cap S_j = \emptyset$ for all $i \neq j$ and $S_1 \cup S_2 \cup \cdots \cup S_k = S$.
	
	\subsection{Elements of graph theory}
	We denote by $\calG = (V,E)$ a directed graph with vertex set $V = \nset{n}$ and edge set $E \subseteq V \times V$.
	We define the graph $\calG = (V,E)$ to be {\em undirected} if $(i,j) \in E$ whenever $(j,i) \in E$. 
	In that case, the order of $i$ and $j$ does not matter, and we interpret the edge set $E$ as the set of unordered pairs $\{i,j\}$, where $(i,j) \in E$.
	Moreover, an undirected graph $\calG = (V,E)$ is called a {\em bipartite graph} if there exist nonempty disjoint subsets $X$ and $Y$ of $V$ such that $X \cup Y = V$ and $\{i,j\} \in E$ only  if $i \in X$ and $j \in Y$.
	Such a bipartite graph is denoted by $G = (X,Y,E_{XY})$, where we denote the edge set by $E_{XY}$ to stress that it contains edges $\{i,j\}$ with $i \in X$ and $j \in Y$. 
	In this paper, we will use the symbol $\calG$ for directed graphs and $G$ for bipartite graphs.
	A set of $t$ edges $m \subseteq E_{XY}$ is called a {\em $t$-matching} in $G$, if no two distinct edges in $m$ share the same vertex. 
	In the special case that $|X| = |Y| = t$, such a  $t$-matching is called a {\em perfect matching.}

	\subsection{Pattern matrices and structured systems}
	By a pattern matrix, we mean a matrix with entries in the set of symbols $\{0,\ast,?\}$. 
	The set of all $p \times q$ pattern matrices will be denoted by $\{0,\ast,?\}^{p \times q}$.
	For a given $p \times q$ pattern matrix $\mathcal{M}$, we define the \emph{pattern class} of $\mathcal{M}$ as
	\begin{equation*}
	\begin{aligned}
	\mathcal{P}(\mathcal{M}) :=  \{ M \in \mathbb{R}^{p \times q} \mid  & M_{ij} = 0 \text{ if } \mathcal{M}_{ij} = 0, \\& M_{ij} \neq 0 \text{ if } \mathcal{M}_{ij} = \ast \}.
	\end{aligned}
	\end{equation*}
	This means that for a matrix $M \in \mathcal{P}(\mathcal{M})$, the entry $M_{ij}$ is either (i) \emph{zero} if $\mathcal{M}_{ij} = 0$, (ii) \emph{nonzero} if $\mathcal{M}_{ij} = \ast$, or (iii) \emph{arbitrary} (zero or nonzero) if $\mathcal{M}_{ij} = \: ?$.
	
	Let $\mathcal{A} \in \{0,\ast,?\}^{n \times n}$ and $\mathcal{B} \in \{0,\ast,?\}^{n \times m}$ be two pattern matrices.
	Consider the linear dynamical system
	\begin{equation}\label{eq:gss}		
	\dot{x}(t) = A x(t) + B u(t),
	\end{equation}
	where  $x \in \mathbb{R}^n$ is the state, and $u \in \mathbb{R}^m$ is the input.
	We will call the family of systems \eqref{eq:gss} with $A \in \calP(\calA)$ and $B \in \calP(\calB)$ a {\em structured system}.
	We denote this structured system by the ordered pair of pattern matrices $(\calA,\calB)$, and we denote by $(A,B)$ a particular system of the form \eqref{eq:gss}.
	Thus, \[(\calA,\calB) = \set{(A,B)}{\bbm A & B\ebm \in \calP(\bbm \calA& \calB\ebm)}.\]
	\section{Problem formulation}\label{s:3}
	
	In this section, we will introduce the problem to be considered in this paper.
	Let $(\calA,\calB)$ be the structured system associated  with $\calA \in \{0,\ast,?\}^{n \times n}$ and $\calB \in \{0,\ast,?\}^{n \times m}$.
	The structured system $(\calA,\calB)$ is called strongly structurally controllable if each $(A,B)$ in this family is controllable.
	In \cite{JHHK2019} necessary and sufficient conditions for strong structural controllability were established.
	Note that in the set-up of \cite{JHHK2019}, all $\ast$ and $?$-entries in $\calA$ and $\calB$ take their values independently. 
	In the present paper we will extend the results of \cite{JHHK2019} and impose constraints on the $\ast$ and $?$-entries.
	In particular, instead of considering the entire family  $(\calA,\calB)$, we will zoom in on a subclass of $(\calA,\calB)$ containing those systems $(A,B)$ that satisfy the condition that a prior given entries in $\bbm A & B \ebm$ are equal.
	We will now formalize this equality constraints on the $\ast$ and $?$-entries.
	
	To do so, consider a pattern matrix $\calM \in \{0,\ast,?\}^{p \times q}$.
	Define the sets of locations of $\ast$ and $?$-entries in  $\calM$ as 
	\[
	\calI_\calM(\ast) := \{(i,j) \in \{1, 2, \ldots, p\} \times \{1,2, \ldots,q\} \mid \calM_{ij} = \ast \}
	\]
	and
	\[
	\calI_\calM(?) := \{(i,j) \in \{1, 2, \ldots, p\} \times \{1,2, \ldots,q\} \mid \calM_{ij} = ? \}.
	\]
	Let 
	$
	\pi^\ast := \{\calI^\ast_{1}, \calI^\ast_{2},\ldots,\calI^\ast_{k}\}
	$
	and 
$
	\pi^? := \{\calI^?_{1}, \calI^?_{2},\ldots,\calI^?_{l}\}
	$
	be partitions of $\calI_{\calM} (\ast)$ and $\calI_{\calM} (?)$, respectively.
	We then call $\pi := \pi^\ast \cup \pi^?$  a {\em coloring} of the pattern matrix $\calM $ and the pair $(\calM, \pi)$ a {\em colored pattern matrix}. 
	Next, we define the so-called {\em colored pattern class} associated with $(\calM, \pi)$ as 
	\[
	\bali
	\calP(\calM, \pi) := \{ &M \in \calP(\calM) \mid  M_{ij} = M_{kl} \mbox{ if there exists }  \\  &\text{$r \in \nset{k}
		$ such that } (i,j), (k,l) \in \calI^\ast_r \\
	\text{ or } & s \in \nset{l} \text{ such that } (i,j), (k,l) \in  \calI^?_s \}.
	\eali
	\]
	In order to visualize the coloring $\pi$, two $\ast$-entries in the same subset $\calI^\ast_r$ are said to have the same color, which will be denoted by the symbol $c_r$.
	Likewise, two $?$-entries in the subset $\calI_{s}^{?}$ are said to have the same color, and this color will be denoted by the symbol $g_s$.
	In this paper, we will always use symbols $c_r~ ( r = 1, 2, \ldots, k)$ for colors associated with $\ast$-entries,  and $g_s~ ( s = 1, 2, \ldots, l)$ for colors associated with $?$-entries.
	With a slight abuse of notation, sometimes we will also use the symbols $c_i$ and $g_i$ to denote nonzero or arbitrary real variables.
		
	\bex \label{ex:pm}
	Consider the colored pattern matrix $(\calM,\pi)$ 
	with
	\beq \label{eq:ex}
	\calM = \resizebox{0.35\hsize}{!}{ $\bbm
		0    & 0    & \ast  & 0    & 0 & \ast &  0 \\
		0    & ?    &0      &\ast  &?  &\ast  &\ast \\
		\ast & 0    &?      &0     &0  & 0    &0\\
		? & ? &\ast   &\ast  &0  & 0    &0\\
		\ast &\ast  &0      &0     &0  & 0    &0\\
		\ebm$} 
	\text{~,~~~}
	\pi = \{\calI^\ast_1, \calI^\ast_2,\calI^?_{1},\calI^?_{2}\},
	\eeq
	where
	\beqn
	\bali
	\calI^\ast_1 &= \{(1,3), (1,6), (2,7), (3,1), (4,3), (4,4)\},\\
	\calI^\ast_2 &= \{(2,4), (2,6),  (5,1), (5,2)\},\\
	\calI^?_{1}  &= \{(4,1),(4,2),(2,5)\},\\
	\calI^?_{2}  &= \{(2,2),(3,3)\}.
	\eali
	\eeqn
	In this example, the $\ast$-entries with locations in $\calI^\ast_1$ have color $c_1$, and those with locations in $\calI^\ast_2$ have color $c_2$.
	Besides, the $?$-entries with locations in $\calI^?_1$ have color $g_1$ and those with locations in $\calI^?_2$ have color $g_2$. 
	
	Thus, $(\calM,\pi)$ can be  visualized by
	\beqn 
	\resizebox{0.45\hsize}{!}{$    \bbm    0 & 0 & \mathcolor{myc1}{c_1} & 0 & 0 & \mathcolor{myc1}{c_1} &  0 \\
		0 & \mathcolor{myg2}{g_2} &0 &\mathcolor{myc2}{c_2} &\mathcolor{myg1}{g_1} &\mathcolor{myc2}{c_2}  &\mathcolor{myc1}{c_1} \\
		\mathcolor{myc1}{c_1}& 0 &\mathcolor{myg2}{g_2} &0 &0& 0 &0\\
		\mathcolor{myg1}{g_1}& \mathcolor{myg1}{g_1}&\mathcolor{myc1}{c_1}&\mathcolor{myc1}{c_1}&0& 0&0\\
		\mathcolor{myc2}{c_2}&\mathcolor{myc2}{c_2}&0&0&0&0&0\\
		\ebm$.}
	\eeqn 
	The corresponding colored pattern class consists of all real matrices of the form
	\beq \label{eq:matrix}
	\resizebox{0.45\hsize}{!}{
		$\bbm
		0          & 0            & c_1       & 0     & 0         & c_1 &0 \\
		0          & g_2    & 0         & c_2   & g_1 &c_2  &c_1 \\
		c_1        & 0            &g_2  & 0     &0          & 0   &0\\
		g_1  &g_1     & c_1       &c_1    &0          & 0   &0\\
		c_2        & c_2          &0          & 0     &0          & 0   &0\\
		\ebm$},
	\eeq
	where the $g_i$ are arbitrary real numbers, and the $c_i$ are arbitrary nonzero real numbers.
	\eex

	A colored pattern matrix $(\calM, \pi)$ is said to have {\em full row rank} if every matrix $M \in \calP(\calM, \pi)$ has full row rank.
	
	Consider now the colored pattern matrix $( \resizebox{0.12\hsize}{!}{$\bbm \calA & \calB \ebm$}, \pi)$ associated with the pattern matrices 
	$\calA \in \{0,\ast,?\}^{n \times n} \mbox{ and } \calB \in \{0,\ast,?\}^{n \times m}$
	and the coloring $\pi = \{ \calI^\ast_1, \calI^\ast_2, \ldots, \calI^\ast_k, \calI^?_1, \calI^?_2, \ldots, \calI^?_l\}$.
	We then define the {\em colored structured system} associated with $( \resizebox{0.12\hsize}{!}{$\bbm \calA & \calB \ebm$}, \pi)$ as 
	\beqn
	(\calA,\calB,\pi) := \set{(A,B)}{\resizebox{0.12\hsize}{!}{$\bbm A & B\ebm$} \in \calP(\resizebox{0.12\hsize}{!}{$\bbm \calA& \calB\ebm$}, \pi)}.
	\eeqn
	We say that this colored structured system is  {\em strongly structurally controllable} if $(A,B)$ is controllable for all $\resizebox{0.12\hsize}{!}{$\bbm A & B \ebm$} \in \calP(\resizebox{0.12\hsize}{!}{$\bbm\calA & \calB \ebm$}, \pi)$.
	We will then simply say that $(\calA,\calB,\pi)$ is {\em controllable}.
	For example, the colored structured system $(\calA,\calB,\pi)$ with $\resizebox{0.12\hsize}{!}{$\bbm \calA & \calB \ebm$}$ and $\pi$ given by \eqref{eq:ex} is controllable, as will be shown later on in this paper.
	The problem that we will investigate in this paper can now be stated as follows.
	\bprb
	Given a colored structured system $(\calA, \calB, \pi)$, find conditions under which it is controllable.
	\eprb 
	\bre \label{r:compa}
    There is a strong relation between the work in this paper and that in \cite{JTBC20181} and \cite{JTBC2018} on controllability of systems on colored graphs. Stated in terms of pattern matrices, the work in \cite{JTBC20181} and \cite{JTBC2018} deals with the very special case of colored structured systems $(\calA,\calB,\pi)$ in which
    \begin{enumerate}
    	\item all diagonal entries of $\calA$ are $?$,
    	\item all off-diagonal entries of $\calA$ are $\ast$ or $0$,
    	\item in $\calB$, each column contains exactly one $\ast$ and each row contains at most one $\ast$,
    	\item the coloring $\pi^? = \{(1,1), (2,2), \ldots, (n,n)\}$ of the $?$-entries is given, i.e.,    	
    	the $?$-entries have distinct colors.
    \end{enumerate}
    In \cite{JTBC20181} and \cite{JTBC2018}, conditions were obtained for controllability of this special class of systems. In the present paper these results will be generalized to general colored structured systems.
	\ere

	\section{Algebraic conditions for controllability}\label{s:4}
	
	In this section, we will provide a sufficient algebraic condition for controllability.
The condition states that a colored structured system is controllable if two particular colored pattern matrices associated with this system have full row rank.

	Let $(\calA, \calB, \pi)$ be a colored structured system with $\calA \in \{0,\ast,?\}^{n \times n}, \calB \in \{0,\ast,?\}^{n \times m}$, and $$\pi = \{ \calI^\ast_1, \calI^\ast_2, \ldots, \calI^\ast_k, \calI^?_1, \calI^?_2, \ldots, \calI^?_l\}.$$ 
	In order to state our first result, for the given $(\calA,\calB,\pi)$ we define an associated new colored pattern matrix $(\bbm \bar{\calA} & \calB\ebm, \bar{ \pi})$ as follows.
	\bdfn \label{d:coloredpatternmatrix}
	We define $\bar{\calA}$ to  be the pattern matrix obtained from $\calA$ by modifying the diagonal entries of $\calA$ as follows
	\[
	\bar{\calA}_{ii} := 
	\begin{cases}
	\ast & \text{ if $\calA_{ii} = 0$,}\\
	?     & \text{ otherwise.}
	\end{cases}
	\]
	The matrix $\calB$ remains unchanged.
	Next, for $r = 1,2, \ldots, k$ and $s = 1,2, \ldots, l$ we remove the diagonal locations from $\calI_{r}^{\ast}$ and $\calI_{s}^{?}$ 
by defining
 $$\bar{\calI}^\ast_r := \{(i,j) \in \calI^\ast_r \mid   i \neq j \}$$
    and
     $$\bar{\calI}^?_s := \{(i,j) \in \calI^?_s \mid  i \neq j \}.$$
	Note that some of the $\bar{\calI}^\ast_r$ or $\bar{\calI}^?_s$ might be empty.
	Next, we partition the set of diagonal locations into $n$ subsets.
	More explicitly, if $i_1, i_2, \ldots, i_w \in \nset{n}$ are the indices such that $\bar{\calA}_{i_j i_j} = \ast$, then for $j = 1, 2, \ldots, w$ we define $$\bar{\calI}^\ast_{k+j} := \{(i_j,i_j)\}.$$ 
	Furthermore, if $t_1, t_2, \ldots, t_{n-w} \in \nset{n}$ are the indices such that $\bar{\calA}_{t_j t_j} = ?$
	for $j = 1, 2, \ldots, n-w$, then we define $$\bar{\calI}^?_{l+j} := \{(t_j,t_j)\} \text{ for $j = 1, 2, \ldots, n-w$}.$$
	Thus we obtain a partition 
	\[\bar{\pi} := \{\bar{\calI}^\ast_1, \bar{\calI}^\ast_2, \ldots, \bar{\calI}^\ast_{k+w}, \bar{\calI}^?_1, \bar{\calI}^?_2, \ldots, \bar{\calI}^?_{l+n-w}\}\]
	of the sets $\calI_{\resizebox{0.08\hsize}{!}{$\bbm \bar{\calA} & \calB \ebm$}}(\ast)$ and $\calI_{\resizebox{0.08\hsize}{!}{$\bbm \bar{\calA} & \calB \ebm$}}(?)$ and define this to be the new coloring $\bar{\pi}$.
	\edfn
	
	We are now ready to state our first result.
	\bthm \label{t:AT}
	Let  $(\calA, \calB, \pi)$ be a colored structured system, and let $(\resizebox{0.12\hsize}{!}{$\bbm \bar{\calA} & \calB \ebm$}, \bar{\pi})$ be the colored pattern matrix obtained from $( \resizebox{0.12\hsize}{!}{ $\bbm \calA & \calB \ebm$}, \pi)$ as in Definition \ref{d:coloredpatternmatrix}.
	Then, $(\calA, \calB, \pi)$ is controllable if both $( \resizebox{0.12\hsize}{!}{$\bbm \calA & \calB \ebm$}, \pi)$ and $(\resizebox{0.12\hsize}{!}{$\bbm \bar{\calA} & \calB \ebm$}, \bar{\pi})$ have full row rank.
	\ethm
	\BP
	The proof of this theorem can be given by slightly adapting that of the sufficient condition in \cite[Theorem 8]{JHHK2019} and is hence omitted.
	\EP
	We will now illustrate Theorem \ref{t:AT} by an example.
	\bex \label{ex:coloredssytem}
	Consider $(\calA,\calB,\pi)$ with $\bbm \calA & \calB \ebm$ and $\pi$ given by \eqref{eq:ex} in Example \ref{ex:pm}.
	Using Definition \ref{d:coloredpatternmatrix}, we obtain the colored pattern matrix $(\resizebox{0.12\hsize}{!}{$\bbm \bar{\calA} & \calB \ebm$}, \bar{\pi})$ with 
	\beq \label{eq:ex1}
	\bali
	\resizebox{0.12\hsize}{!}{ $\bbm \bar{\calA} & \calB \ebm$}  &= 
	\resizebox{0.35\hsize}{!}{ $\bbm
		\ast & 0    & \ast  & 0    & 0    & \ast &0 \\
		0    & ?    &0      &\ast  &?     &\ast  &\ast \\
		\ast & 0    &?      &0     &0     & 0    &0 \\
		?    & ?    &\ast   &?     &0     & 0    &0 \\
		\ast &\ast  &0      &0     &\ast  & 0    &0 \\
		\ebm$}  ,\\  
	\bar{\pi} &= \{\bar{\calI}^\ast_1, \bar{\calI}^\ast_2, \bar{\calI}^\ast_3, \bar{\calI}^\ast_4,  \bar{\calI}^?_{1},\calI^?_{2} ,\calI^?_{3} ,\calI^?_{4} ,\calI^?_{5}\},
	\eali
	\eeq
	where 
	\beqn
	\bali
	\bar{\calI}^\ast_1 &= \{(1,3), (1,6), (2,7), (3,1), (4,3)\},\\
	\bar{\calI}^\ast_2 &= \{ (2,4), (2,6),  (5,1), (5,2)\},
	\bar{\calI}^\ast_3 = \{(1,1)\} , 
	\bar{\calI}^\ast_4 = \{(5,5)\},\\
	\bar{\calI}^?_{1}  &= \{(4,1),(4,2),(2,5)\}, 
	\bar{\calI}^?_{2}  = \emptyset, 
	\bar{\calI}^?_{3}  = \{(2,2)\}, \\ 
	\bar{\calI}^?_{4}  &= \{(3,3)\},  
	\bar{\calI}^?_{5}  = \{(4,4)\}.
	\eali
	\eeqn	
	It turns out that both $( \resizebox{0.12\hsize}{!}{$\bbm \calA & \calB \ebm$}, \pi)$ and $(\resizebox{0.12\hsize}{!}{$\bbm \bar{\calA} & \calB \ebm$}, \bar{\pi})$ have full row rank.
	Indeed, let $M$ be a matrix in $\calP(\resizebox{0.12\hsize}{!}{$\bbm \calA & \calB\ebm $}, \pi)$.
	Then, $M$ is of the form \eqref{eq:matrix}.
	Let $M'$ be the submatrix of $M$ obtained by removing the fourth and fifth column from $M$.
	It is easily seen that $\det(M') = -c_1^4 c_2$, which is nonzero for $c_1 \text{  and } c_2$.
	Hence, all matrices $M$ given by \eqref{d:coloredpatternmatrix} have full row rank, so all $(\resizebox{0.12\hsize}{!}{$\resizebox{0.12\hsize}{!}{$\bbm \calA & \calB\ebm $}$}, \pi)$ has full row rank.
	Similarly, one can verify that $(\resizebox{0.12\hsize}{!}{$\bbm \bar{\calA} & \calB \ebm$}, \bar{\pi})$ has full row rank, so, by Theorem \ref{t:AT}, we conclude that $(\calA,\calB,\pi)$ is controllable.
	\eex

\bre
In \cite{JHHK2019}, necessary and sufficient conditions for controllability of structured systems $(\calA,\calB)$ without a coloring were established, also in terms of two rank tests.
We note that Theorem \ref{t:AT} generalizes this result to colored systems. 
The conditions obtained in the present paper are however only sufficient.
To illustrate that the conditions in Theorem \ref{t:AT} are not necessary for controllability, we will provide a counterexample of a colored structured system that is controllable while one of the conditions does not hold.
\ere	
	\bex
	Consider $(\mathcal{A},\mathcal{B},\pi)$ with
	$\resizebox{0.12\hsize}{!}{$\bbm \calA &\calB \ebm$} = \bbm  \ast & \ast& \ast \\ \ast & 0 & \ast\ebm$
	and 
	$\pi = \{\calI^\ast_1, \calI^\ast_2\},$
	where $\calI^\ast_1 = \{(1,1), (1,2), (2,1)\} \text{ and } \calI^\ast_2 = \{(1,3), (2,3)\}.$	
	The corresponding colored pattern class consists of all matrices of the form
	$
	\bbm  c_{1} & c_{1}& c_{2}\\ c_{1} & 0 & c_{2}\ebm,
	$
	where $c_{1},c_{2}$ are nonzero real numbers.
	The matrix $\bbm B & AB \ebm$ is equal to $	\bbm c_{2}& 2c_{1}c_{2}\\ c_{2}& c_{1}c_{2} \ebm$ which has full row rank for every choice of $c_1$ and $c_2$.
	By the Kalman rank test, we conclude that $(A,B)$ is controllable for all $\bbm A & B\ebm \in \calP(\resizebox{0.12\hsize}{!}{$\bbm \calA &\calB \ebm$}, \pi)$, i.e., $(\mathcal{A},\mathcal{B},\pi)$ is controllable.
	
	Next, we will show that the second condition in Theorem \ref{t:AT} is not satisfied.
	Let $(\resizebox{0.12\hsize}{!}{$\bbm \bar{\calA} &\calB \ebm$}, \bar{\pi})$ be the colored pattern matrix obtained from $(\resizebox{0.12\hsize}{!}{$\bbm \calA &\calB \ebm$}, \pi)$ as in Definition \ref{d:coloredpatternmatrix} with
	\[
	\resizebox{0.12\hsize}{!}{$\bbm \bar{\calA} &\calB \ebm$} = \bbm  ? & \ast& \ast \\ \ast & \ast & \ast\ebm
	,
	\bar{\pi} = \{\bar{\calI}^\ast_1, \bar{\calI}^\ast_2, \bar{\calI}^\ast_3, \bar{\calI}^?_1\},
	\]
	where $\bar{\calI}^\ast_1 = \{(1,2), (2,1)\}$, $\bar{\calI}^\ast_2 = \{(1,3), (2,3)\}$, $\bar{\calI}^\ast_3 = \{(2,2)\}$ and $\bar{\calI}^?_1 = \{(1,1)\}$.
	Now, consider the matrix
	$M = \begin{bmatrix}1&1&1\\ 1&1&1\end{bmatrix}.$
	Clearly, $M \in \calP(\resizebox{0.12\hsize}{!}{$\bbm \bar{\calA} &\calB \ebm$}, \bar{\pi})$ while it does not have full row rank.
	Hence, we conclude that the second condition in Theorem \ref{t:AT} is not satisfied.
	\eex

	Checking whether a colored pattern matrix has full row rank is in general not an easy task.
	Therefore, in the sequel we will develop a test for this in terms of a so-called color change rule on the graph associated with the colored pattern matrix.
	In order to do this, in the next section, we will consider {\em square} colored pattern matrices and establish graph-theoretic conditions under which all matrices in the associated pattern class are nonsingular.
	
	\section{Conditions for nonsingularity of square colored pattern matrices} \label{s:5}
	Let $\calN \in \{0,\ast,?\}^{t \times t}$ be a square pattern matrix.
	We define the {\em pattern class} of $\calN$ as
	\beqn
	\bali
	\calP(\calN) := \{N \in \bbC^{t \times t} \mid  N_{ij} = 0 &\mbox{ if } \calN_{ij} = 0, \\  N_{ij} \neq 0 &\mbox{ if } \calN_{ij} = \ast\}.\\
	\eali
	\eeqn
	Note that here and in the sequel, in the context of pattern classes for square pattern matrices, we will allow complex matrices.
	Let $\pi = \{\calI^\ast_1, \calI^\ast_2, \ldots, \calI^\ast_k, \calI^?_1, \calI^?_2, \ldots, \calI^?_l\}$ be a coloring of $\calN$.
	Again, two $\ast$-entries in the same subset $I^\ast_r$ are said to have the same color, visualized by a symbol $c_r$, and two $?$-entries in the same subset $I^?_s$ are said to have the same color, visualized by a symbol $g_s$.
	As before, $(\calN, \pi)$ denotes the colored pattern matrix associated with $\calN$ and $\pi$.
	The corresponding {\em pattern class} of $(\calN, \pi)$ is given by
	\beqn
	\bali
	\calP(\calN, \pi) = &\{N \in \calP(\calN) \mid  N_{ij} = N_{mn} \mbox{ if there exists} \\ & r \in \nset{k} \text{ such that } (i,j), (m,n) \in \calI^\ast_r \\ \text{ or } & s \in \nset{l} \text{ such that } (i,j), (m,n) \in \calI^?_s\}.
	\eali
	\eeqn
	We say that $(\calN, \pi)$ is {\em nonsingular} if all matrices in $\calP(\calN, \pi)$ are nonsingular.
	In this section, we will establish necessary and sufficient conditions for  nonsingularity in terms of bipartite graphs.
	
	We define the bipartite graph $G = (X,Y,E_{XY})$ associated with the $t \times t$ pattern matrix $\calN$ as follows.
	Take as vertex sets $X = \{x_1, x_2, \ldots, x_t\}$ and $Y = \{y_1, y_2, \ldots, y_t\}$.
	An edge $\{x_i, y_j\}$ belongs to the edge set $E_{XY}$  if $\calN_{ji} = \ast$ or $?$.
	To distinguish the edges corresponding to entries equal to $?$ and $\ast$, we introduce the two subsets 
	\[
	E^\ast_{XY} := \{ \{x_i, y_j\} \in E_{XY} \mid \calN_{ji} = \ast\}\]
	and
	\[ E^?_{XY} := \{ \{x_i, y_j\} \in E_{XY} \mid \calN_{ji} = ?\}.
	\]
	To visualize these different kinds of edges, we use solid lines to represent the edges in $E^\ast_{XY}$ and dashed lines to represent the edges in $E^?_{XY}$.
	In addition, the coloring $\pi$ induces a partition of the edge set $E_{XY}$:
	$$ \pi_{XY} := \{ E^{\ast 1}_{XY}, E^{\ast 2}_{XY}, \ldots, E^{\ast k}_{XY}, E^{? 1}_{XY}, E^{? 2}_{XY}, \ldots, E^{? l}_{XY}\}$$  
	in which for $r = 1,2, \ldots, k$
	\[
	E^{\ast r}_{XY} := \{ \{x_i,y_j\} \in E^\ast_{XY} \mid (j, i) \in \calI^\ast_r \},
	\]
	and for $s = 1,2, \ldots, l$
	\[
	E^{? s}_{XY} := \{ \{x_i,y_j\} \in E^?_{XY} \mid (j, i) \in \calI^?_s \}.
	\]
    The partition $\pi_{XY}$ is a coloring of the edge set $E_{XY}$.
    The edges in the same subset $E^{\ast r}_{XY}$ inherit the color $c_r$ corresponding to $\calI^\ast_{r}$.
    Similarly, the edges in the subset $E^{? s}_{XY}$ inherit the color $g_s$ corresponding to $\calI^?_{s}$.
	Thus we define the {\em colored bipartite graph} associated with $(\calN ,\pi)$ as $G(\calN, \pi) =(X,Y,E_{XY},\pi_{XY})$.
	
	\bex \label{ex:CBG}
	Consider the square colored pattern matrix $(\calN, \pi)$ with 
	\[\calN = \bbm
	\ast& 0 &? \\
	?& ?&\ast\\
	\ast&\ast&0\\
	\ebm
	,~~
	\pi = \{\calI^\ast_{1},\calI^\ast_{2},\calI^?_{1}, \calI^?_{2}\},\]
	where 
	$\calI^\ast_{1} = \{(1,1), (2,3)\}$, 
	$\calI^\ast_{2} = \{(3,1), (3,2)\}$,
	$\calI^?_{1} = \{(2,1),$ $(2,2)\}$ and $\calI^?_{2} = \{(1,3)\}$. 
	The associated colored bipartite graph $G(\calN, \pi)$ is depicted in Figure \ref{f:CBG}.
	\begin{figure}[!ht]
		\centering
		\scalebox{0.8}{\begin{tikzpicture}[scale=0.9]
			\tikzset{VertexStyle1/.style = {shape = circle,
					color=black,
					fill=white!93!black,
					minimum size=0.5cm,
					text = black,
					inner sep = 2pt,
					outer sep = 1pt,
					minimum size = 0.55cm}
			}	
			\tikzset{VertexStyle2/.style = {shape = circle,
					color=black,
					fill=black!93!white,
					minimum size=0.5cm,
					text = white,
					inner sep = 2pt,
					outer sep = 1pt,
					minimum size = 0.55cm}
			}	
			\tikzset{VertexStyle3/.style = {shape = circle,
					color=white,
					fill=white,
					minimum size=0.5cm,
					text = black,
					inner sep = 2pt,
					outer sep = 1pt,
					minimum size = 0.55cm}
			}	
			\node[VertexStyle1,draw](1) at (0,1.5) {$1$};
			\node[VertexStyle1,draw](2) at (0,0) {$2$};
			\node[VertexStyle1,draw](3) at (0,-1.5) {$3$};
			\node[VertexStyle1,draw](4) at (2,1.5) {$1$};
			\node[VertexStyle1,draw](5) at (2,0) {$2$};
			\node[VertexStyle1,draw](6) at (2,-1.5) {$3$};
			\node[VertexStyle3,draw](7) at (0,2.25) {$X$};
			\node[VertexStyle3,draw](8) at (2,2.25) {$Y$};
			\Edge[ style = {-,> = latex',pos=0.5},color = myc1, label = $c_{1}$,labelstyle={inner sep=0pt}](1)(4);
			\Edge[ style = {-,> = latex',pos=0.2},color = myc1, label = $c_{1}$,labelstyle={inner sep=0pt}](3)(5);
			\Edge[ style = {-,> = latex',pos=0.2},color = myc2, label = $c_{2}$,labelstyle={inner sep=0pt}](1)(6);			
			\Edge[ style = {-,> = latex',pos=0.7},color = myc2, label = $c_{2}$,labelstyle={inner sep=0pt}](2)(6);
			\Edge[ style = {-,> = latex',pos=0.5,dashed},color = myg1, label = $g_{1}$,labelstyle={inner sep=0pt}](1)(5);
			\Edge[ style = {-,> = latex',pos=0.2,dashed},color = myg1, label = $g_{1}$,labelstyle={inner sep=0pt}](2)(5);	
			\Edge[ style = {-,> = latex',pos=0.2,dashed},color = myg2, label = $g_{2}$,labelstyle={inner sep=0pt}](3)(4);			
			\end{tikzpicture}}
		\caption{Colored bipartite graph $G(\calN, \pi)$.}
		\label{f:CBG}
	\end{figure}
	\eex

	In order to proceed, we will now review some concepts associated with perfect matchings in bipartite graphs. (See also \cite{JTBC20181,JTBC2018}.)
	Let $p$ be a perfect matching in $G(\calN, \pi)$.  
	The \emph{spectrum} of $p$ is defined as the set of colors of the edges in $p$.
	More explicitly, if the perfect matching $p$ is given by 
	\beq \label{eq:PM}
	p =\{\{x_1,y_{\gamma(1)}\}, \{x_2,y_{\gamma(2)}\},\ldots,\{x_t,y_{\gamma(t)}\} \},
	\eeq 
	where $\gamma$ denotes a permutation of $(1, 2, \ldots,t)$, and $$c_{i_1}, c_{i_2}, \ldots, c_{i_j},g_{i_{j+1}}, g_{i_{j+2}}, \ldots, g_{i_{t}}  \text{ with } j \leq t $$ are the respective colors of the edges in $p$, then the spectrum of $p$ is defined as $\{c_{i_1}, c_{i_2}, \ldots, c_{i_j}, g_{i_{j+1}}, g_{i_{j+2}}, \ldots, g_{i_{t}}\}$, where the same color can appear multiple times.	
	We say that two perfect matchings are \emph{equivalent} if they have the same spectrum. 
	Obviously, this leads to a partition of the set of all perfect matchings of $G(\calN, \pi)$ into {\em equivalence classes}. 
	We denote these equivalence classes of perfect matchings by $\bbP_{1}, \bbP_{2}, \dots, \bbP_{r}$ in which perfect matchings in the same class $\bbP_i$ are equivalent.
	Define the {\em spectrum of the equivalence class $\bbP_i$} to be the (common) spectrum of the perfect matchings in this class and denote it by $\spec(\bbP_i)$.
	Clearly, for $i \neq j$, we have $\bbP_i \cap \bbP_j = \emptyset$ and $\spec(\bbP_i) \neq \spec(\bbP_j)$.
	The \emph{sign} of the perfect matching $p$ given by \eqref{eq:PM} is defined as $\text{sign}(p)=(-1)^{m}$, where $m$ is the number of swaps required to permute $(1,2, \ldots,t)$ to $(\gamma(1),\gamma(2),\ldots,\gamma(t))$. 
	Finally, we define the signature of $\mathbb{P}_{i}$ as  $\text{sgn}(\mathbb{P}_{i}) :=\sum_{p\in\mathbb{P}_{i}}\text{sign}(p).$
	In other words, the signature of the equivalence class $\bbP_i$ is equal to the sum of the signs of the perfect matchings contained in $\bbP_i$.
	In order to illustrate the above definitions, we now give an example.
	
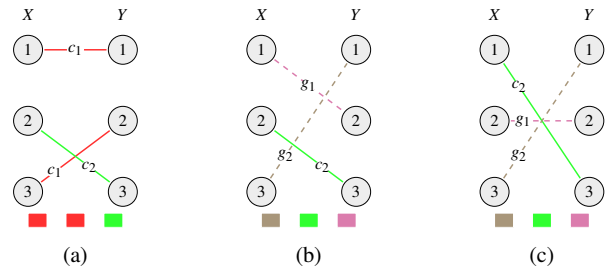
\begin{figure}[!ht]
		\begin{subfigure}{.13\textwidth}\centering
			\scalebox{0.7}{\begin{tikzpicture}[scale=0.9]
				\tikzset{VertexStyle1/.style = {shape = circle,
						color=black,
						fill=white!93!black,
						minimum size=0.5cm,
						text = black,
						inner sep = 2pt,
						outer sep = 1pt,
						minimum size = 0.55cm}
				}	
				\tikzset{VertexStyle2/.style = {shape = circle,
						color=black,
						fill=black!93!white,
						minimum size=0.5cm,
						text = white,
						inner sep = 2pt,
						outer sep = 1pt,
						minimum size = 0.55cm}
				}
				\tikzset{VertexStyle3/.style = {shape = circle,
						color=white,
						fill=white,
						minimum size=0.5cm,
						text = black,
						inner sep = 2pt,
						outer sep = 1pt,
						minimum size = 0.55cm}
				}	
				\node[VertexStyle1,draw](1) at (0,1.5) {$1$};
				\node[VertexStyle1,draw](2) at (0,0) {$2$};
				\node[VertexStyle1,draw](3) at (0,-1.5) {$3$};
				\node[VertexStyle1,draw](4) at (2,1.5) {$1$};
				\node[VertexStyle1,draw](5) at (2,0) {$2$};
				\node[VertexStyle1,draw](6) at (2,-1.5) {$3$};
				\node[VertexStyle3,draw](7) at (0,2.25) {$X$};
				\node[VertexStyle3,draw](8) at (2,2.25) {$Y$};
				\Edge[ style = {-,> = latex',pos=0.5},color = myc1, label = $c_{1}$,labelstyle={inner sep=0pt}](1)(4);
				\Edge[ style = {-,> = latex',pos=0.2},color = myc1, label = $c_{1}$,labelstyle={inner sep=0pt}](3)(5);
				\Edge[ style = {-,> = latex',pos=0.7},color = myc2, label = $c_{2}$,labelstyle={inner sep=0pt}](2)(6);
				\fill[fill=myc1,draw=white] (0,-2.25) rectangle +(4mm,3mm);
				\fill[fill=myc1,draw=white] (0.8,-2.25) rectangle +(4mm,3mm);
				\fill[fill=myc2,draw=white] (1.6,-2.25) rectangle +(4mm,3mm);
				\end{tikzpicture}}
			\caption{}
		\end{subfigure}
		\qquad
		\begin{subfigure}{.13\textwidth}\centering
			\scalebox{0.7}{\begin{tikzpicture}[scale=0.9]
				\tikzset{VertexStyle1/.style = {shape = circle,
						color=black,
						fill=white!93!black,
						minimum size=0.5cm,
						text = black,
						inner sep = 2pt,
						outer sep = 1pt,
						minimum size = 0.55cm}
				}	
				\tikzset{VertexStyle2/.style = {shape = circle,
						color=black,
						fill=black!93!white,
						minimum size=0.5cm,
						text = white,
						inner sep = 2pt,
						outer sep = 1pt,
						minimum size = 0.55cm}
				}	
				\tikzset{VertexStyle3/.style = {shape = circle,
						color=white,
						fill=white,
						minimum size=0.5cm,
						text = black,
						inner sep = 2pt,
						outer sep = 1pt,
						minimum size = 0.55cm}
				}	
				\node[VertexStyle1,draw](1) at (0,1.5) {$1$};
				\node[VertexStyle1,draw](2) at (0,0) {$2$};
				\node[VertexStyle1,draw](3) at (0,-1.5) {$3$};
				\node[VertexStyle1,draw](4) at (2,1.5) {$1$};
				\node[VertexStyle1,draw](5) at (2,0) {$2$};
				\node[VertexStyle1,draw](6) at (2,-1.5) {$3$};
				\node[VertexStyle3,draw](7) at (0,2.25) {$X$};
				\node[VertexStyle3,draw](8) at (2,2.25) {$Y$};
				\Edge[ style = {-,> = latex',pos=0.7},color = myc2, label = $c_{2}$,labelstyle={inner sep=0pt}](2)(6);
				\Edge[ style = {-,> = latex',pos=0.5,dashed},color = myg1, label = $g_{1}$,labelstyle={inner sep=0pt}](1)(5);
				\Edge[ style = {-,> = latex',pos=0.2,dashed},color = myg2, label = $g_{2}$,labelstyle={inner sep=0pt}](3)(4);
				\fill[fill=myg2,draw=white] (0,-2.25) rectangle +(4mm,3mm);
				\fill[fill=myc2,draw=white] (0.8,-2.25) rectangle +(4mm,3mm);
				\fill[fill=myg1,draw=white] (1.6,-2.25) rectangle +(4mm,3mm);
				\end{tikzpicture}}
			\caption{ }
		\end{subfigure}
		\qquad
		\begin{subfigure}{.13\textwidth}\centering
			\scalebox{0.7}{\begin{tikzpicture}[scale=0.9]
				\tikzset{VertexStyle1/.style = {shape = circle,
						color=black,
						fill=white!93!black,
						minimum size=0.5cm,
						text = black,
						inner sep = 2pt,
						outer sep = 1pt,
						minimum size = 0.55cm}
				}	
				\tikzset{VertexStyle2/.style = {shape = circle,
						color=black,
						fill=black!93!white,
						minimum size=0.5cm,
						text = white,
						inner sep = 2pt,
						outer sep = 1pt,
						minimum size = 0.55cm}
				}	
				\tikzset{VertexStyle3/.style = {shape = circle,
						color=white,
						fill=white,
						minimum size=0.5cm,
						text = black,
						inner sep = 2pt,
						outer sep = 1pt,
						minimum size = 0.55cm}
				}	
				\node[VertexStyle1,draw](1) at (0,1.5) {$1$};
				\node[VertexStyle1,draw](2) at (0,0) {$2$};
				\node[VertexStyle1,draw](3) at (0,-1.5) {$3$};
				\node[VertexStyle1,draw](4) at (2,1.5) {$1$};
				\node[VertexStyle1,draw](5) at (2,0) {$2$};
				\node[VertexStyle1,draw](6) at (2,-1.5) {$3$};
				\node[VertexStyle3,draw](7) at (0,2.25) {$X$};
				\node[VertexStyle3,draw](8) at (2,2.25) {$Y$};
				\Edge[ style = {-,> = latex',pos=0.2},color = myc2, label = $c_{2}$,labelstyle={inner sep=0pt}](1)(6);			
				\Edge[ style = {-,> = latex',pos=0.2,dashed},color = myg1, label = $g_{1}$,labelstyle={inner sep=0pt}](2)(5);	
				\Edge[ style = {-,> = latex',pos=0.2,dashed},color = myg2, label = $g_{2}$,labelstyle={inner sep=0pt}](3)(4);
				\fill[fill=myg2,draw=white] (0,-2.25) rectangle +(4mm,3mm);
				\fill[fill=myc2,draw=white] (0.8,-2.25) rectangle +(4mm,3mm);
				\fill[fill=myg1,draw=white] (1.6,-2.25) rectangle +(4mm,3mm);
				\end{tikzpicture}}
			\caption{ }
		\end{subfigure}
		\caption{ (a) Perfect matching $p_{1}$ with spectrum  $\{ c_1,c_1, c_2\}$ and $\text{sign}(p_{1})=-1$. (b) Perfect matching $p_{2}$ with spectrum $ \{ g_2,c_2, g_1\}$ and $\text{sign}(p_{2})=1$. (c) Perfect matching $p_{3}$ with spectrum $\{ g_2,c_2, g_1\}$ and $\text{sign}(p_{3})=-1$.}
		\label{f:SCP}
	\end{figure}

	\bex \label{ex:CG}
	Consider the colored bipartite graph $G(\calN, \pi)=(X,Y,E_{XY},\pi_{XY})$ depicted in Figure \ref{f:CBG}.
	It has three perfect matchings $p_{1}$, $p_{2}$ and $p_{3}$ in $G(\calN, \pi)$, which are depicted in Figures \ref{f:SCP}(a), \ref{f:SCP}(b) and  \ref{f:SCP}(c), respectively. Clearly, $p_{2}$ and $p_{3}$ are equivalent. Thus, the equivalence classes are $\mathbb{P}_{1}=\{p_{1} \}$ and $\bbP_{2}=\{p_{2},p_{3}\}$ with signature $\text{sgn}(\mathbb{P}_{2})=\text{sign}(p_{2})+\text{sign}(p_{3})=0$ and $\text{sgn}(\mathbb{P}_{1})=-1$.
	\eex
	We are now ready to state a necessary and sufficient condition for a square colored pattern matrix to be nonsingular.
	\bthm\label{t:GTSM}
	Let $(\calN, \pi)$ be a square colored pattern matrix and $G(\calN, \pi) = (X, Y, E_{XY}, \pi_{XY})$ its associated bipartite graph. Then, $(\calN, \pi)$ is nonsingular if and only if in $G(\calN, \pi)$
	the following three conditions hold:
	\begin{enumerate}
		\item there exists at least one perfect matching,
		\item exactly one equivalence class of perfect matchings has a nonzero signature, 
		\item the spectrum of this equivalence class contains only colors corresponding to edges in $E^\ast_{XY}$, i.e., solid edges.
	\end{enumerate}
	\ethm
	The proof can be found in the Appendix.
	Theorem \ref{t:GTSM} is a generalization of \cite[Theorem 8]{JTBC2018} which provides a necessary and sufficient condition for the special case of square colored pattern matrices only containing $0$ and $\ast$-entries.
	\bex
	Consider the square colored pattern matrix $(\calN, \pi)$ given in Example \ref{ex:CBG} and its associated graph $G(\calN, \pi)$ depicted in Figure \ref{f:CBG}.
	In Example \ref{ex:CG}, it has already been shown that $G(\calN, \pi)$ admits exactly one equivalence class, $\bbP_1 = \{p_1\}$, with nonzero signature. Moreover, $\spec(\bbP_1) =\{c_1, c_1, c_2\}$, which only contains colors associated with solid edges.
	Therefore, by Theorem \ref{t:GTSM}, we conclude that $(\calN, \pi)$ is nonsingular.
	\eex

	\section{Color change rule and graph-theoretic conditions}\label{s:6}
	In this section, we will establish a graph-theoretic test for checking whether a colored pattern matrix has full row rank. 
	This test will be in terms of a so-called color change rule on the associated graph.
	Color change rules for checking the rank of a pattern matrix have been studied before, see e.g., \cite{AIM2008},\cite{MZC2014},\cite{TD2015},\cite{JHHK2019}, \cite{JTBC2018}.
	Here, we will start off with introducing a new color change rule tailored for our purpose.
	
	Let $(\calM, \pi)$ be the colored pattern matrix associated with $\calM \in \{0,\ast,?\}^{p \times q}$ ($p \leq q$) and $\pi = \{\calI^\ast_1, \calI^\ast_2, \ldots,\calI^\ast_k, \calI^?_1, \calI^?_2, \ldots, \calI^?_l\}$.
	Define a directed graph $\calG(\calM, \pi)=(V,E)$ associated with $(\calM, \pi)$ as follows.
	Take the vertex set $V$ equal to $\{1, 2, \ldots,q \}$. 
	Define the edge set $E\subseteq V\times V$ as
	\[
	E := \set{(i,j)}{\calM_{ji} = \ast \text{ or } \calM_{ji} = ?}.
	\]
	The coloring $\pi$ gives the following partition of the edge set $E$:
	$$\pi_E := \{E^\ast_1, E^\ast_2, \ldots, E^\ast_k, E^?_{1}, E^?_{2}, \ldots, E^?_{l}\}$$ 
	in which for $r = 1,2, \ldots, k$
	\[E^\ast_{r}: = \{ (i,j) \in E \mid (j,i) \in \calI^\ast_{r}\},\]
	and for $s = 1,2, \ldots, l$
	\[E^?_{s}: = \{ (i,j) \in E \mid (j,i) \in \calI^?_{s}\}.\]
	We call the partition $\pi_E$ a coloring of the edge set $E$.
	To visualize the coloring $\pi_E$,  for $r = 1,2, \ldots, k$ we represent the edges in $E^\ast_r$ by solid arrows with color $c_r$ (inherited from $\calI_{r}^\ast$).
	For $s = 1 , 2, \ldots, l$ we represent the edges in $E^?_{s}$ by dashed arrows with color $g_s$  (inherited from $\calI^?_s$).
	Thus, we obtain a colored graph $\calG(\calM, \pi) = (V,E,\pi_E)$ associated with $(\calM, \pi)$.
	Colored graphs were studied before in \cite{JTBC2018}.
	In order to illustrate the above, we provide an example.
	\bex
	Consider the colored pattern matrix $(\calM, \pi)$ of Example \ref{ex:pm}.
	The associated graph $\calG(\calM, \pi)$ is depicted in Figure \ref{f:GCP}.
	\begin{figure}[!ht]
		\centering
		\scalebox{0.85}{\begin{tikzpicture}[scale=0.5]
			\tikzset{VertexStyle1/.style = {shape = circle,
					color=black,
					fill=white!93!black,
					minimum size=0.5cm,
					text = black,
					inner sep = 2pt,
					outer sep = 1pt,
					minimum size = 0.55cm}
			}	
			\tikzset{VertexStyle2/.style = {shape = circle,
					color=black,
					fill=black!93!white,
					minimum size=0.5cm,
					text = white,
					inner sep = 2pt,
					outer sep = 1pt,
					minimum size = 0.55cm}
			}	
			\node[VertexStyle1,draw](1) at (-2,2){\textbf{1}};         
			\node[VertexStyle1,draw](2) at(-2,-2){\textbf{2}};					
			\node[VertexStyle1,draw](3) at (4,-2){\textbf{3}};
			\node[VertexStyle1,draw](4) at (4,2){\textbf{4}};
			\node[VertexStyle1,draw](5) at (1.5,4.5){\textbf{5}};			
			\node[VertexStyle1,draw](6) at (-5,2){\textbf{6}};
			\node[VertexStyle1,draw](7) at (-5,-2){\textbf{7}};
			\Loop[ style = {->,> = latex',pos=0.5, out = -60, in = -120, dashed},color=myg2, label = $g_{2}$,labelstyle={inner sep=0pt},dist=1.3cm](2);
			\Loop[ style = {->,> = latex',pos=0.5, out = -60, in = -120, dashed},color=myg2, label = $g_{2}$,labelstyle={inner sep=0pt},dist=1.3cm](3);
			\Loop[ style = {->,> = latex',pos=0.5, out = 30, in =-30},color=myc1, label = $\textcolor{black}{c_{1}}$,labelstyle={inner sep=0pt},dist=1.3cm](4);
			\Edge[style = {->,> = latex',pos = 0.8},color=myc1, label = $c_{1}$,labelstyle={inner sep=0pt}](1)(3);
			\Edge[style = {->,> = latex',pos = 0.5, bend right}, color=myc1, label = $c_1$, labelstyle={inner sep=0pt}](3)(1);
			\Edge[style = {->,> = latex',pos = 0.5}, color=myc1, label = $c_1$, labelstyle={inner sep=0pt}](6)(1);
			\Edge[style = {->,> = latex',pos = 0.5}, color=myc1, label = $c_1$, labelstyle={inner sep=0pt}](7)(2);
			\Edge[style = {->,> = latex',pos = 0.5}, color=myc1, label = $c_1$, labelstyle={inner sep=0pt}](3)(4);
			\Edge[style = {->,> = latex',pos = 0.3}, color=myc2, label = $c_2$, labelstyle={inner sep=0pt}](1)(5);
			\Edge[style = {->,> = latex',pos = 0.5, out = -120, in = 0}, color=myc2, label = $c_2$, labelstyle={inner sep=0pt}](4)(2);
			\Edge[style = {->,> = latex',pos = 0.5}, color=myc2, label = $c_2$, labelstyle={inner sep=0pt}](6)(2);
			\Edge[style = {->,> = latex',pos = 0.2, out = 40, in = -90}, color=myc2, label = $c_2$, labelstyle={inner sep=0pt}](2)(5);
			\Edge[style = {->,> = latex',pos = 0.8, dashed}, color=myg1, label = $g_1$, labelstyle={inner sep=0pt}](2)(4);
			\Edge[style = {->,> = latex',pos = 0.7, dashed}, color=myg1, label = $g_1$, labelstyle={inner sep=0pt}](1)(4);
			\Edge[style = {->,> = latex',pos = 0.2,dashed}, color=myg1, label = $g_1$, labelstyle={inner sep=0pt}](5)(2);		
			\end{tikzpicture}}
		\caption{The colored graph $\calG(\calM, \pi)$ associated with $(\calM, \pi)$.}
		\label{f:GCP}
	\end{figure}
	\eex
	
	We will now introduce a color change rule for $\calG(\calM,\pi)$.
	In this graph, initially all vertices are colored ``white.''
	The color change rule will prescribe under what conditions vertices will change their color to ``black.''
	In earlier work, color change rules usually deal with conditions under which a {\em single} vertex colors a {\em single} white neighboring vertex black. (See \cite{MZC2014}, \cite{TD2015}, \cite{JHHK2019} and the references therein.)
	In the present paper we will deal with {\em sets} of vertices that color {\em sets} of vertices black. (See also \cite{JTBC2018}.)
	We will now describe this rule.
	Let $X$ and $Y$ be two nonempty subsets of the vertex set $V$, containing the same number of vertices, i.e., $|X| = |Y|$.
	Define an associated colored bipartite graph $G(\pi) = (X,Y,E_{XY},\pi_{XY})$ as follows:
	\[
	E_{XY} := \set{\{x_i,y_j\}}{ x_i \in X, y_j \in Y, (x_i,y_j) \in E}.
	\]
	Obviously, the partition $\pi_E$ of $E$ induces a partition $$\pi_{XY} = \{E^{\ast 1}_{XY}, E^{\ast 2}_{XY}, \ldots, E^{\ast k}_{XY},  E^{? 1}_{XY}, E^{? 2}_{XY}, \ldots, E^{? l}_{XY}\}$$ of $E_{XY}$ by defining 
	for $r = 1,2, \ldots, k$
	$$
	E^{\ast r}_{XY} := \{\{x_i,y_j\} \in E_{XY} \mid (x_i,y_j) \in E^\ast_r\},
    $$
	and
	for $s = 1,2, \ldots, l$
	$$
	E^{? s}_{XY} := \{\{x_i,y_j\} \in E_{XY} \mid (x_i,y_j) \in E^?_s\}.
	$$
	Note that  some of these sets might be empty.
	Removing all the empty sets, we then obtain a partition
	\[
	\pi_{XY} = \{E^{\ast i_1}_{XY}, E^{\ast i_2}_{XY}, \ldots, E^{\ast i_w}_{XY},  E^{? j_1}_{XY}, E^{? j_2}_{XY}, \ldots, E^{? j_v}_{XY}\}
	\]
	of $E_{XY}$ with $w \leq k \text{ and } v \leq l$.
	The edges in $E^{\ast i_r}_{XY}$ have color $c_{i_r}$, and the edges in $E^{? j_r}_{XY}$ have color $g_{j_{r}}$.
	Without loss of generality, we renumber $c_{i_1}, c_{i_2}, \ldots, c_{i_w}$ as $c_{1}, c_{2}, \ldots, c_{w}$ and $g_{j_{1}}, g_{j_{2}}, \ldots, g_{j_{v}}$ as $g_1, g_2, \ldots, g_{v}$.
	
	Next, return to the colored graph $\calG(\calM, \pi) = (V,E)$.
	Suppose that all vertices in $V$ are colored either black or white.
	Take two nonempty subsets $X$ and $Y$ of the vertex set $V$.
	We say that  $Y$ is a \emph{color-perfect white neighbor} of $X$ if:
	\begin{enumerate}
		\item $Y$ and $X$ contain the same number of vertices, i.e., $|Y| = |X|$.
		\item $Y$ is equal to the set of white out-neighbors of $X$, i.e., 
		$$ Y = \{y_j \in V \mid y_j \mbox{ is white and } (x_i,y_j) \in E \mbox{ for some } x_i \in X \};$$
		\item in the associated bipartite graph $G(X,Y,E_{XY},\pi_{XY})$,
		there exists a perfect matching, exactly one equivalence class of perfect matchings has a nonzero signature, and the spectrum of this equivalence class only contains colors corresponding to edges in $E^\ast_{XY}$, i.e., solid edges in $G(X,Y,E_{XY},\pi_{XY})$.
	\end{enumerate} 
	Based on the notion of color-perfect white neighbor, we now introduce a {\em color change rule} as follows. 
	\begin{enumerate}
		\item Initially, all vertices in $V$ are colored white.
		\item If there exist two vertex sets $Y\subseteq \{1, 2, \ldots, p\}$ and $X \subseteq \{1,2,\ldots, q\}$ such that $Y$ is a color-perfect white neighbor of $X$, then change the color of all vertices in $Y$ to black.
		\item Repeat the second step until no further color changes are possible.
	\end{enumerate}
    We define a {\em derived set $\calD$} as a set of black vertices in $V$ obtained by following the procedure above.
    Note that derived sets are not unique but may depend on the subsequent choices of $Y$ and $X$ in the second step above.
    An illustrative example can be found in \cite[Example 26]{JTBC2018}.
    The graph $\calG(\calM, \pi)$ is called \emph{colorable} if there exists a derived set $\calD$ such that $\calD = \{1,2, \ldots, p\}$.
	Of course, the remaining vertices $\{p+1, p+2, \ldots,q \}$ can never be colored black, since they have no incoming edges.

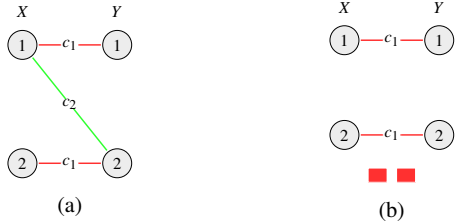
\begin{figure}[!ht]
		\begin{subfigure}{0.23\textwidth}\centering
			\scalebox{0.7}{\begin{tikzpicture}[scale=0.9]
				\tikzset{VertexStyle1/.style = {shape = circle,
						color=black,
						fill=white!93!black,
						minimum size=0.5cm,
						text = black,
						inner sep = 2pt,
						outer sep = 1pt,
						minimum size = 0.55cm}
				}	
				\tikzset{VertexStyle2/.style = {shape = circle,
						color=black,
						fill=black!93!white,
						minimum size=0.5cm,
						text = white,
						inner sep = 2pt,
						outer sep = 1pt,
						minimum size = 0.55cm}
				}	
				\tikzset{VertexStyle3/.style = {shape = circle,
						color=white,
						fill=white,
						minimum size=0.5cm,
						text = black,
						inner sep = 2pt,
						outer sep = 1pt,
						minimum size = 0.55cm}
				}	
				\node[VertexStyle1,draw](1) at (0,2) {$1$};
				\node[VertexStyle1,draw](2) at (0,-0.5) {$2$};
				\node[VertexStyle1,draw](3) at (2,2) {$1$};
				\node[VertexStyle1,draw](4) at (2,-0.5) {$2$};
				\node[VertexStyle3,draw](5) at (0,2.7) {$X$};
				\node[VertexStyle3,draw](6) at (2,2.7) {$Y$};
				\Edge[ style = {-, >= latex',pos=0.5},color = myc1, label = $c_{1}$,labelstyle={inner sep=0pt}](1)(3);
				\Edge[ style = {-,> = latex',pos=0.5},color = myc2, label = $c_{2}$,labelstyle={inner sep=0pt}](1)(4);
				\Edge[ style = {-,> = latex',pos=0.5},color = myc1, label = $c_{1}$,labelstyle={inner sep=0pt}](2)(4);
				\end{tikzpicture}}
			\caption{}
		\end{subfigure}
		\begin{subfigure}{0.23\textwidth}\centering
			\scalebox{0.7}{\begin{tikzpicture}[scale=0.9]
				\tikzset{VertexStyle1/.style = {shape = circle,
						color=black,
						fill=white!93!black,
						minimum size=0.5cm,
						text = black,
						inner sep = 2pt,
						outer sep = 1pt,
						minimum size = 0.55cm}
				}	
				\tikzset{VertexStyle2/.style = {shape = circle,
						color=black,
						fill=black!93!white,
						minimum size=0.5cm,
						text = white,
						inner sep = 2pt,
						outer sep = 1pt,
						minimum size = 0.55cm}
				}	
				\tikzset{VertexStyle3/.style = {shape = circle,
						color=white,
						fill=white,
						minimum size=0.5cm,
						text = black,
						inner sep = 2pt,
						outer sep = 1pt,
						minimum size = 0.55cm}
				}
				\node[VertexStyle1,draw](1) at (0,2) {$1$};
				\node[VertexStyle1,draw](2) at (0,0) {$2$};
				\node[VertexStyle1,draw](3) at (2,2) {$1$};
				\node[VertexStyle1,draw](4) at (2,0) {$2$};
				\node[VertexStyle3,draw](5) at (0,2.7) {$X$};
				\node[VertexStyle3,draw](6) at (2,2.7) {$Y$};
						\Edge[ style = {-, >= latex',pos=0.5},color = myc1, label = $c_{1}$,labelstyle={inner sep=0pt}](1)(3);
				\Edge[ style = {-,> = latex',pos=0.5},color = myc1, label = $c_{1}$,labelstyle={inner sep=0pt}](2)(4);
				\fill[fill=myc1,draw=white] (0.5,-1) rectangle +(4mm,3mm);
				\fill[fill=myc1,draw=white] (1.1,-1) rectangle +(4mm,3mm);
				\end{tikzpicture}}
			\caption{}
		\end{subfigure}
		\caption{(a) Colored bipartite graph $G(X,Y,E_{XY},\pi_{XY})$ associated with $X = \{6,7\}$ and $Y=\{1,2\}$. (b) Perfect matching $p_{1}$ with spectrum $\{ c_1, c_1\}$ and $\text{sign}(p_{1})= 1$.}
		\label{f:XY}
	\end{figure}

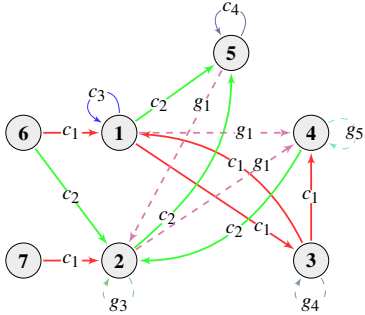
\begin{figure}[!ht]
		\centering
		\scalebox{0.85}{\begin{tikzpicture}[scale=0.5]
			\tikzset{VertexStyle1/.style = {shape = circle,
					color=black,
					fill=white!93!black,
					minimum size=0.5cm,
					text = black,
					inner sep = 2pt,
					outer sep = 1pt,
					minimum size = 0.55cm}
			}	
			\tikzset{VertexStyle2/.style = {shape = circle,
					color=black,
					fill=black!93!white,
					minimum size=0.5cm,
					text = white,
					inner sep = 2pt,
					outer sep = 1pt,
					minimum size = 0.55cm}
			}	
			\node[VertexStyle1,draw](1) at (-2,2){\textbf{1}};         
			\node[VertexStyle1,draw](2) at(-2,-2){\textbf{2}};					
			\node[VertexStyle1,draw](3) at (4,-2){\textbf{3}};
			\node[VertexStyle1,draw](4) at (4,2){\textbf{4}};
			\node[VertexStyle1,draw](5) at (1.5,4.5){\textbf{5}};			
			\node[VertexStyle1,draw](6) at (-5,2){\textbf{6}};
			\node[VertexStyle1,draw](7) at (-5,-2){\textbf{7}};
			\Loop[ style = {->,> = latex',pos=0.5, out = 90, in = 150},color=myc3, label = $\textcolor{black}{c_{3}}$,labelstyle={inner sep=0pt},dist=1.3cm](1);
			\Loop[ style = {->,> = latex',pos=0.5, out = -60, in = -120, dashed},color=myg3, label = $\textcolor{black}{g_{3}}$,labelstyle={inner sep=0pt},dist=1.3cm](2);
			\Loop[ style = {->,> = latex',pos=0.5, out = -60, in = -120, dashed},color=myg4, label = $\textcolor{black}{g_{4}}$,labelstyle={inner sep=0pt},dist=1.3cm](3);
			\Loop[ style = {->,> = latex',pos=0.5, out = 30,  in =-30,   dashed},color=myg5, label = $\textcolor{black}{g_{5}}$,labelstyle={inner sep=0pt},dist=1.3cm](4);	
			\Loop[ style = {->,> = latex',pos=0.5, out = 60, in = 120},color=myc4, label = $\textcolor{black}{c_{4}}$,labelstyle={inner sep=0pt},dist=1.3cm](5);
			\Edge[style = {->,> = latex',pos = 0.8},color=myc1, label = $c_{1}$,labelstyle={inner sep=0pt}](1)(3);
			\Edge[style = {->,> = latex',pos = 0.5, bend right}, color=myc1, label = $c_1$, labelstyle={inner sep=0pt}](3)(1);
			\Edge[style = {->,> = latex',pos = 0.5}, color=myc1, label = $c_1$, labelstyle={inner sep=0pt}](6)(1);
			\Edge[style = {->,> = latex',pos = 0.5}, color=myc1, label = $c_1$, labelstyle={inner sep=0pt}](7)(2);
			\Edge[style = {->,> = latex',pos = 0.5}, color=myc1, label = $c_1$, labelstyle={inner sep=0pt}](3)(4);
			\Edge[style = {->,> = latex',pos = 0.3}, color=myc2, label = $c_2$, labelstyle={inner sep=0pt}](1)(5);
			\Edge[style = {->,> = latex',pos = 0.5, out = -120, in = 0}, color=myc2, label = $c_2$, labelstyle={inner sep=0pt}](4)(2);
			\Edge[style = {->,> = latex',pos = 0.5}, color=myc2, label = $c_2$, labelstyle={inner sep=0pt}](6)(2);
			\Edge[style = {->,> = latex',pos = 0.2, out = 40, in = -90}, color=myc2, label = $c_2$, labelstyle={inner sep=0pt}](2)(5);
			\Edge[style = {->,> = latex',pos = 0.8, dashed}, color=myg1, label = $g_1$, labelstyle={inner sep=0pt}](2)(4);
			\Edge[style = {->,> = latex',pos = 0.7, dashed}, color=myg1, label = $g_1$, labelstyle={inner sep=0pt}](1)(4);
			\Edge[style = {->,> = latex',pos = 0.2,dashed}, color=myg1, label = $g_1$, labelstyle={inner sep=0pt}](5)(2);
			\end{tikzpicture}}
		\caption{The colored graph $\calG([ \bar{\calA} ~~ \calB ], \bar{\pi})$.}
		\label{f:GCP2}
	\end{figure}

	\bex \label{ex:MColo}
	Consider $(\calM, \pi)$ given by \eqref{eq:ex} and its associated graph $\calG(\calM, \pi)$ depicted in Figure \ref{f:GCP}.
	Initially, color all vertices white.
	First, let $X = \{6,7\}$ and $Y = \{1,2\}$.
	It turns out that  $Y$ is a color-perfect white neighbor of $X$.
	Indeed, in the associated colored bipartite graph $G = (X,Y,E_{XY},\pi_{XY})$ depicted in Figure \ref{f:XY}, there exists exactly one equivalence class $\bbP_1 = \{p_1\}$ with nonzero signature and  $\spec(\bbP_1) = \{c_1, c_1\}$.
	Consequently, we change the color of vertices $1$ and $2$ to black.
	Next, let $X'  = \{1,2,3\}$ and $Y'  = \{3,4,5\}$. 
	Then $Y' $ is a color-perfect white neighbor of $X'$.
	Indeed, the associated colored bipartite graph $G=(X',Y',E_{X'Y'},\pi_{X'Y'})$ is depicted in Figure \ref{f:CBG}.
	In Example \ref{ex:CG}, we have shown that in the bipartite graph $G=(X',Y',E_{X'Y'},\pi_{X'Y'})$ there exists exactly one equivalence class of perfect matchings with nonzero signature and the spectrum of this equivalence class contains only colors corresponding to solid edges.
	Therefore, vertices $3,4,5$ are colored black, and we conclude that $\calG(\calM, \pi)$ is colorable.
	\eex
	We now arrive at the main result of this section.
	\bthm \label{t:GTM}
	Let $(\calM, \pi)$ be the colored pattern matrix associated with  $\calM \in \{0,\ast,?\}^{p \times q}$ ($p \leq q$)  and  $$\pi = \{\calI^\ast_1, \calI^\ast_2, \ldots,\calI^\ast_k, \calI^?_1, \calI^?_2, \ldots,\calI^?_l\}.$$
	Then, $(\calM, \pi)$ has full row rank if its associated graph $\calG(\calM, \pi)$ is colorable.
	\ethm
	In order to prove this theorem, we need the following instrumental result.
	\blem \label{l:colorchange}
	Let $(\calM, \pi)$ be a colored pattern matrix and $\calG(\calM,\pi)$ its associated colored graph.
	Suppose that the vertices $1, 2, \ldots, p$ are  black or white, and those in $p+1, p+2, \ldots, q$ are  white.
	Define the diagonal matrix $D \in\mathbb{R}^{p\times p}$ by
	\beqn\label{eq:D5}
	D_{ii} := \begin{cases}
		1& \text{if the vertex $i$ is black} ,\\
		0& \text{otherwise} .
	\end{cases}
	\eeqn
	Suppose further that $Y =\{y_{1}, y_{2},\ldots,y_{r} \} \subseteq\{1,2,\ldots,p \}$ is a color-perfect white neighbor of $X = \{x_{1}, x_{2},\ldots,x_{r} \} \subseteq \{1,2, \ldots, q\}$.
	Define the diagonal matrix $\Delta \in\mathbb{R}^{p\times p}$ by
	\[
	\Delta:= \sum_{i=1}^{r} e_{y_i}e_{y_i}^T ,
	\]  
	where $e_{y_i}$ denotes the $y_i$th column of the $p \times p$ identity matrix $I$.
	Then for every $M \in \calP(\calM, \pi)$ we have that  $\begin{bmatrix} M & D\end{bmatrix}$ has full row rank if and only if $\begin{bmatrix} M & D+ \Delta \end{bmatrix}$ has full row rank.
	\elem
	\BP
	The `only if' part is trivial and is hence omitted.
	To prove the `if' part, suppose that, for all $M \in \calP(\calM, \pi)$ the matrix $\begin{bmatrix} M & D+ \Delta\end{bmatrix}$ has full row rank.
	Let $z \in \bbR^p$ be such that $z^T \bbm M & D\ebm = 0$.
	In the sequel, for a given vector $z \in \bbR^p$ and a given index set $\alpha = \{\alpha_1, \alpha_2, \ldots, \alpha_r\} \subseteq \nset{p}$, we define the vector 
	$$z_{\alpha} := (z_{\alpha_1}, z_{\alpha_2}, \ldots, z_{\alpha_r})^T.$$
	Analogously, for a given matrix $M \in \bbR^{p \times q}$ and two given index sets $\alpha = \{\alpha_1, \alpha_2, \ldots, \alpha_r\} \subseteq \nset{p}$ and $\beta = \{\beta_1, \beta_2, \ldots, \beta_s\} \subseteq \nset{q}$, we define the matrix $M_{\alpha \beta}$ by
	$
	(M_{\alpha \beta} )_{ij} := M_{\alpha_i \beta_j}.
	$
	In what follows, we aim to show that $z_Y = \0$.
	Indeed, if $z_Y = \0$ then $z^T \bbm M & D + \Delta \ebm = z^T \bbm M & D\ebm = \0,$
	which would prove that $z = \0$.
	So, $\bbm M & D\ebm$ has full row rank.
	To show that, indeed, $z_Y = \0$, let $\alpha$ be the set of black vertices. 
	Clearly, it holds that $\alpha\subseteq \{1,2,\ldots,p\}$ and $\alpha \cap Y = \emptyset$.
	Since $z^T \bbm M & D \ebm = \0$, we then obtain
	\[
	z_Y^T M_{YX} + z_\alpha^T M_{\alpha X} + z_\beta^T M_{\beta X}  = \0 \text{ and } z_\alpha = \0,
	\]
	where $\beta = \{1,2,\ldots,p\} \setminus (Y \cup \alpha)$.
	Since $Y$ is a color-perfect white neighbor of $X$, by Theorem \ref{t:GTSM} we have that $M_{YX}$ is nonsingular and $M_{\beta X} = \0$.
	This implies that $z_Y$ must be equal to $\0$.
	This completes the proof.
	\EP
	We are now ready to prove Theorem~\ref{t:GTM}.
	\BP[Proof of Theorem \ref{t:GTM}]
	Suppose that $\calG(\calM, \pi)$ is colorable. Let $M \in \calP(\calM, \pi)$. 
	By repeatedly applying Lemma~\ref{l:colorchange}, it follows that $M$ has full row rank if and only if $\bbm M & I \ebm$ has full row rank, which is obviously true.
	Therefore, we conclude that $M$ has full row rank, which completes the proof.
	\EP
	To show that the condition in Theorem~\ref{t:GTM} is not a necessary condition, we provide the following counterexample.
	\bex\label{ex:count}
	Consider the colored pattern matrix $(\calM, \pi)$ with
	\[
	\calM = \begin{bmatrix}
	\ast & \ast & \ast & 0\\
	0   & \ast & 0 & \ast\\
	\ast & 0 & \ast & \ast\\
	\end{bmatrix}, ~~
	\pi = \{\calI^\ast_{1},\calI^\ast_{2},\calI^\ast_{3}\},
	\]
	where $\calI^\ast_{1} = \{(1,1),(3,1)\}$, $\calI^\ast_{2} = \{(1,2), (2,2), (2,4), (3,3)\}$ and $\calI^\ast_{3} = \{(1,3),(3,4)\}$.
	It will turn out that the associated graph  $\calG(\calM, \pi)$ depicted in Figure \ref{f:t7} is not colorable, while $(\calM, \pi)$ has full row rank.
	Indeed, one can verify that none of the subsets of $\{1,2,3,4\}$ has a color-perfect white neighbor. Hence, the graph $\calG(\calM, \pi)$ is not colorable. 
	However, all matrices of the form
	\[\resizebox{0.28\hsize}{!}{ $\bbm
	c_{1} & c_{2} & c_{3} & 0\\
	0     & c_{2} & 0     & c_{2}\\
	c_{1} & 0     & c_{2} & c_{3}\\
	\ebm$}
	\]
	with $c_i \neq 0$ for $i = 1,2,3$ have full row rank.
	Indeed, by taking
	$$P = \resizebox{0.18\hsize}{!}{ $\bbm 1 & 0 & 0\\ 0 & 1 & 0\\ -1 & 1 & 1\ebm$} \text{ and } Q = \resizebox{0.18\hsize}{!}{ $\bbm 1 & 0 & 0 & 0 \\ 0 & 1 & 0 & 0\\ 0 & 0 & 1 & 0\\ 0 & 0 & 1 & 1\ebm$}$$ 
	we obtain
	$$P M Q = \bbm
	c_{1} & c_{2} & c_{3} & 0\\
	0     & c_{2} & 0     & c_{2}\\
	0      & 0     & 2c_{2} & c_2 + c_{3}\\
	\ebm,$$
	which clearly has full row rank for all choices of $c_i \neq 0$.
	This provides a counterexample as claimed.
	\begin{figure}[!ht]
		\centering
		\scalebox{0.85}{\begin{tikzpicture}[scale=0.9]
			\tikzset{VertexStyle1/.style = {shape = circle,
					color=black,
					fill=white!96!black,
					minimum size=0.55cm,
					text = black,
					inner sep = 2pt,
					outer sep = 1pt,
					minimum size = 0.55cm}
			}		
			\tikzset{VertexStyle2/.style = {shape = circle,
					color=black,
					fill=black!96!white,
					minimum size=0.55cm,
					text = white,
					inner sep = 2pt,
					outer sep = 1pt,
					minimum size = 0.55cm}
			}
			\tikzset{VertexStyle1/.style = {shape = circle,
					color=black,
					fill=white!96!black,
					minimum size=0.55cm,
					text = black,
					inner sep = 2pt,
					outer sep = 1pt,
					minimum size = 0.55cm}
			}		
			\tikzset{VertexStyle2/.style = {shape = circle,
					color=black,
					fill=black!96!white,
					minimum size=0.55cm,
					text = white,
					inner sep = 2pt,
					outer sep = 1pt,
					minimum size = 0.55cm}
			}			
			\node[VertexStyle1,draw](1) at (-1,0) {$1$};
			\node[VertexStyle1,draw](2) at (1,1.5) {$2$};
			\node[VertexStyle1,draw](3) at (1,-1.5) {$3$};
			\node[VertexStyle1,draw](4) at (3,0) {$4$};
			\Loop[ style = {->,> = latex',pos=0.5,out=150,in= 210},color=myc1, label = $\textcolor{black}{c_{1}}$,labelstyle={inner sep=0pt},dist=1.3cm](1);
			\Loop[ style = {->,> = latex',pos=0.5,out=-60,in=-120},color=myc2, label = $\textcolor{black}{c_{2}}$,labelstyle={inner sep=0pt},dist=1.3cm](3);
			\Loop[ style = {->,> = latex',pos=0.5,out=60 ,in= 120},color=myc2, label = $\textcolor{black}{c_{2}}$,labelstyle={inner sep=0pt},dist=1.3cm](2);
			\Edge[style = {->,> = latex',pos = 0.5, out = 0, in = 120},color=myc1, label = $c_{1}$,labelstyle={inner sep=0pt}](1)(3);
			\Edge[style = {->,> = latex',pos = 0.5, out = -180, in = -60},color=myc3, label = $c_{3}$,labelstyle={inner sep=0pt}](3)(1);
			\Edge[style = {->,> = latex',pos = 0.5},color=myc2, label = $c_{2}$,labelstyle={inner sep=0pt}](2)(1);				
			\Edge[style = {->,> = latex',pos = 0.5},color=myc2, label = $c_{2}$,labelstyle={inner sep=0pt}](4)(2);
			\Edge[style = {->,> = latex',pos = 0.5},color=myc3, label = $c_{3}$,labelstyle={inner sep=0pt}](4)(3);
			\end{tikzpicture}}
		\caption{Colored bipartite graph $\calG(\calM, \pi)$.}
		\label{f:t7}
	\end{figure}
	\eex
	Finally, based on the rank tests in Theorem \ref{t:AT} and the result in Theorem \ref{t:GTM}, we obtain the following sufficient graph-theoretic condition for controllability of colored structured systems.
	\bthm \label{t:GT}
	Consider the colored structured system $(\calA, \calB, \pi)$. 
	Let $(\resizebox{0.12\hsize}{!}{$\bbm \bar{\calA} & \calB \ebm$}, \bar{\pi})$ be the colored pattern matrix associated with $(\calA, \calB, \pi)$ given by Definition \ref{d:coloredpatternmatrix}.
	Then, $(\calA, \calB, \pi)$ is controllable if both 
	graphs $\calG( \resizebox{0.12\hsize}{!}{$\bbm \calA & \calB \ebm$}, \pi)$ and $\calG(\resizebox{0.12\hsize}{!}{$\bbm \bar{\calA} & \calB \ebm$}, \bar{\pi})$ are colorable.
	\ethm
	To conclude this section, we illustrate the above theorem by an example.
	\bex
	Consider the colored structured system $(\calA,\calB,\pi)$ given in Example \ref{ex:coloredssytem}.
	Denote by $\calG(\resizebox{0.12\hsize}{!}{ $ \bbm \calA & \calB\ebm $},\pi )$ and $\calG(\resizebox{0.12\hsize}{!}{$\bbm \bar{\calA} & \calB \ebm$}, \bar{\pi})$ the colored graphs associated with $(\resizebox{0.12\hsize}{!}{ $ \bbm \calA & \calB\ebm $},\pi) $ and $(\resizebox{0.12\hsize}{!}{ $\bbm \bar{\calA} & \calB \ebm$},\bar{ \pi})$.
	In Example \ref{ex:MColo}, we have already shown that $\calG(\resizebox{0.12\hsize}{!}{ $ \bbm \calA & \calB\ebm $}, \pi)$ depicted in Figure \ref{f:GCP} is colorable.
	It remains to show that the graph $\calG(\resizebox{0.12\hsize}{!}{$\bbm \bar{\calA} & \calB \ebm$}, \bar{\pi})$, depicted in Figure \ref{f:GCP2}, is also colorable.
	Clearly, the set $\{1,2\}$ is a color-perfect white neighbor of $\{6,7\}$.
	Hence, we color vertices $1 \text{ and } 2$ black.
	Subsequently, $\{3,4,5\}$ is a color-perfect white neighbor of $\{1,2,3\}$. This means that the vertices $3, 4 \text{ and } 5$ are also colored black.
	Therefore, we find that  $\calG(\resizebox{0.12\hsize}{!}{$\bbm \bar{\calA} & \calB \ebm$}, \bar{\pi})$ is colorable.
	By Theorem \ref{t:GT}, we conclude that $(\calA, \calB, \pi)$ is controllable.
	\eex
	\bre
	Theorem \ref{t:GT} can of course be applied to the special case described in Remark \ref{r:compa}.
	    Indeed, if the colored system $(\calA,\calB,\pi)$ satisfies the conditions $1$ to $4$ in Remark \ref{r:compa}, then it is easily verified that $\bar{\calA} = \calA$ and the new coloring $\bar{ \pi}$ coincides with the original coloring $\pi$.
	    Thus we find that $(\calA,\calB,\pi)$ is controllable if the single colored pattern matrix $(\resizebox{0.12\hsize}{!}{$\bbm \calA & \calB \ebm$},\pi)$ is colorable.
	\ere

\section{Conclusion} \label{S:conclusion}
	In this paper, we have studied strong structural controllability of linear structured systems in which the structure of the system matrices is assumed to be given by zero/nonzero/arbitrary pattern matrices.
	In contrast to the work in \cite{JHHK2019} in which the nonzero and arbitrary entries of the system matrices are completely independent, in the present paper we have studied the general case that certain equality constraints among these entries are given, in the sense that a {\em priori} given entries in the system matrices are restricted to take arbitrary but identical values.
    We have formalized this general class of structured systems by introducing the concepts of colored pattern matrices and colored structured systems.
    In this setup, we have established sufficient algebraic conditions for strong structural controllability of a given colored structured system. These conditions are in terms of a rank test on two colored pattern matrices associated with this system.
    We have shown that these conditions are not necessary by providing an example in which a colored structured system is controllable while the conditions are not satisfied.
    Next, we have developed a graph-theoretic condition for a given colored pattern matrix to have full row rank. This condition involves a new color change rule and the concept of colorability of the graph associated with this pattern matrix.
    To do so, a necessary and sufficient graph-theoretic condition for the nonsingularity of a given square colored pattern matrix has been established.
    Finally, we have established sufficient graph-theoretic conditions  under which a given colored structured system is strongly structurally controllable.   
    
    In this paper, the conditions that we have provided are not necessary, and hence finding necessary {\em and} sufficient conditions is still an open problem.
    In addition, we have focused only on finding conditions for controllability, but providing suitable algorithms  to check these conditions (see, e.g., \cite{WRS2014}) still remains an open problem.
    Finally, other possible future research directions could address system-theoretic concepts like strong targeted controllability \cite{vWCT2017,MCT2015} and identifiability \cite{WTC2018} for colored structured systems.

 \addtolength{\textheight}{-.2875cm}
	
	\appendix
	\section{Proof of Theorem \ref{t:GTSM}} \label{s:appendix}
	\BP
	Denote the dimension of $(\calN, \pi)$ by $t$.
	Let $N\in \calP(\calN, \pi)$. From the well-known Leibniz formula for the determinant, we have
	\beqn
	\det(N)=\sum_{\gamma} \left(\text{sign}(\gamma)\prod_{i=1}^{t}N_{\gamma(i) i }\right),
	\eeqn
	where the sum ranges over all permutations $\gamma$ of $(1, 2, \ldots,t)$, and $\text{sign}(\gamma)=(-1)^{m}$ with $m$ the number of swaps necessary to permute $(1, 2, \ldots,t)$ to $(\gamma(1), \gamma(2),\ldots,\gamma(t))$. 
	Clearly, \[\prod_{i=1}^{t}N_{\gamma(i) i}\neq 0\] only if there exists at least one perfect matching \[p = \{\{1,\gamma(1)\},\{2,\gamma(2)\}, \ldots, \{t,\gamma(t)\}\}\] in $G(\calN, \pi)$. 
	Therefore, we rewrite the Leibniz formula as 
	\beqn \label{eq:detperma}
	\det(N) = \sum_{p} \left(\text{sign}(p)\prod_{i=1}^{t}N_{p(i)i}\right),
	\eeqn
	where $p$ ranges over all perfect matchings in $G(\calN, \pi)$ and  $\text{sign}(p)$ is the sign of the perfect matching $p$ (We now identify perfect matchings with their permutations).
	Suppose that there exist $r$ equivalence classes of perfect matchings $\bbP_{1}, \bbP_{2},\ldots,\bbP_{r}$. 
	Then, we have that 	
	\beq \label{eq:deteqc}
	\det(N) = \sum_{\rho=1}^{r} \left(\text{sgn}(\mathbb{P}_{\rho})\prod_{i=1}^{t}N_{p(i)i} \right),
	\eeq
	where, for $\rho=1, 2, \ldots,r$, in the product appearing in the $\rho$th term, the $p$ is an arbitrary matching in $\mathbb{P}_{\rho}$.
	
	We now move to prove the `if' part. Suppose that there exists at least one perfect matching, exactly one equivalence class of perfect matchings with a nonzero signature, and the spectrum of this equivalence class contains only colors corresponding to solid edges. 
	Without loss of generality, we denote that equivalence class by $\mathbb{P}_{1}$.
	Clearly, for every $N\in\mathcal{P}(\calN, \pi)$, we then have
	\beqn \label{eq:eqcl}
	\det(N) = \text{sgn}(\mathbb{P}_{1})\prod_{i=1}^{t}N_{p(i)i}\neq 0,
	\eeqn
	where $p\in\mathbb{P}_{1}$ is arbitrary. 
	Since the spectrum of $p \in \bbP_1$ only contains colors associated with solid edges (whose symbols correspond  to nonzero values), this implies that $(\calN, \pi)$ is nonsingular. Thus, we have proved the `if' part.
	
	Next, we prove the `only if' part. 
	To do so, suppose that $(\calN, \pi)$ is nonsingular and, for $G(\calN, \pi)$,  at least one of the following statements holds:
	\begin{enumerate} [label=(\roman*)]
		\item There does not exist any perfect matching in $G(\calN, \pi)$.
		\item There does not exist an equivalence class of perfect matchings with a nonzero signature.
		\item There exist at least two equivalence classes of perfect matchings with nonzero signature.	
		\item There exists exactly one equivalence class of perfect matchings with a nonzero signature, but its spectrum contains at least one color corresponding to a dashed edge.	
	\end{enumerate}
	Clearly, both in case (i) and (ii), it is obvious  that $\det N=0$ for all $N \in \calP(\calN, \pi)$. This leads to a contradiction.	
	Consider case (iii). 
	Without loss of generality, suppose $\mathbb{P}_{1}$ and $\mathbb{P}_{2}$ have nonzero signature. 
	The signatures of the remaining equivalence classes can be either zero or nonzero.
	Suppose $c_1, c_2, \ldots, c_k, g_1, g_2, \ldots, g_l$ are the colors associated with the partition $\pi_{XY}$.
	Introduce matrices $N \in \bbC^{t \times t}$ as follows:
	\beqn
	N_{ij}:=\begin{cases}
		a_{r} &\text{ if } (j,i) \mbox{ has color } c_r \text{ for some } r, \\
		a_{k+r} &\text{ if } (j,i)  \mbox{ has color } g_r \text{ for some } r,\\
		0 &\text{ otherwise},
	\end{cases}
	\eeqn
	where  $a_{1}, a_{2}, \ldots, a_{k+ l}$ are independent, nonzero, variables that can take values in $\bbC$.
	Clearly, for all choices of the complex values $a_{1}, a_{2}, \ldots, a_{k+ l}$,  we have $N \in \calP(\calN, \pi)$.
	From formula \eqref{eq:deteqc} for the determinant of $N$, it is clear that the perfect matchings in the equivalence class $\mathbb{P}_{\rho}$ give a contribution 
	$$\text{sgn}(\mathbb{P}_{\rho}) a_{1}^{j_{1}}a_{2}^{j_{2}} \cdots a_{k+l}^{j_{k+l}},$$
	where the degrees corresponding to the multiplicities of the colors of the perfect matchings in  $\mathbb{P}_{\rho}$.	
	By construction, we have $\text{spec}(\mathbb{P}_{1})\neq \text{spec}(\mathbb{P}_{2})$.
	Without loss of generality,  assume that the multiplicity $\epsilon_{1}$ of $c_{1}$ in $\mathbb{P}_{1}$ is different from that in $\mathbb{P}_{2}$, which is denoted by $\epsilon_{2}$.
	Then, $\det(N)$ can be expressed as 
	\beq\label{eq:determ}
	\det(N)= \text{sgn}(\mathbb{P}_{1})\phi_{1}a_{1}^{\epsilon_{1}} + \text{sgn}(\mathbb{P}_{2})\phi_{2}a_{1}^{\epsilon_{2}} + f(a_{1}),
	\eeq
	where $\phi_{1}$ and $\phi_{2}$ are determined by $a_{2},\ldots,a_{k+l}$, and the polynomial $f(a_{1})$ corresponds to the remaining equivalence classes.
	It may happen that some monomials in $f(a_{1})$ contain $a_{1}$ with multiplicity $\epsilon_{1}$ or $\epsilon_{2}$. 
	By taking common factors in \eqref{eq:determ}, we then have
	\beq \label{eq:polynomial}
	\text{sgn}(\mathbb{P}_{1})\psi_{1}a_{1}^{\epsilon_{1}} + \text{sgn}(\mathbb{P}_{2})\psi_{2}a_{1}^{\epsilon_{2}}+f'(a_{1}) = 0,
	\eeq
	where $\psi_{1}$ and $\psi_{2}$ depend on $a_{2},\ldots,a_{k+t}$. In addition, the polynomial $f'(a_{1})$ does not contain the monomials with $a_{1}^{\epsilon_{1}}$ and $a_{1}^{\epsilon_{2}}$. Clearly, the variables $a_{2},\ldots,a_{k+l}$ can be chosen such that $\psi_{1}\neq 0$ and $\psi_{2}\neq 0$. 
	By the fundamental theorem of algebra, we conclude that the polynomial equation \eqref{eq:polynomial} has at least one nonzero complex solution. 
	This implies that for some choice of nonzero complex values $a_{1}, a_{2}, \ldots,a_{k+l}$, we have that $\det N=0$, and hence we reach a contradiction.  	
	
	Finally, consider the case (iv).  
	Suppose that exact one equivalence class of perfect matchings has a nonzero signature, and its spectrum contains at least one color corresponding to some dashed edge. 
	Let  $p = \{\{1,\gamma(1)\}, \{2,\gamma(2)\},\ldots, \{t,\gamma(t)\}\}$ be a perfect matching in $\mathbb{P}_{1}$, where $\gamma$ denotes a permutation on $(1,2, \ldots,t)$. 
	Without loss of generality, assume that the edge $\{1,\gamma(1)\}$ is a dashed edge with color $g_r$ for some $r$.
	The remaining edges in $p$ can be either solid or dashed.
	This implies that  $$\det(N) = \text{sgn}(\mathbb{P}_{1}) \prod_{i = 1}^{t}N_{\gamma(i)i} = \text{sgn}(\bbP_1) \cdot g_r \cdot \prod_{i = 2}^{t}N_{\gamma(i)i},$$
	where $g_r$ represents an arbitrary complex value.
	It is obvious that $\det(N) = 0$ if $g_r$ is chosen as zero.
	Again, we reach a contradiction.
	This completes the proof.	
	\EP
	\bibliographystyle{IEEEtran}
	\bibliography{bibliography}

\begin{thebibliography}{10}
\providecommand{\url}[1]{#1}
\csname url@samestyle\endcsname
\providecommand{\newblock}{\relax}
\providecommand{\bibinfo}[2]{#2}
\providecommand{\BIBentrySTDinterwordspacing}{\spaceskip=0pt\relax}
\providecommand{\BIBentryALTinterwordstretchfactor}{4}
\providecommand{\BIBentryALTinterwordspacing}{\spaceskip=\fontdimen2\font plus
\BIBentryALTinterwordstretchfactor\fontdimen3\font minus
  \fontdimen4\font\relax}
\providecommand{\BIBforeignlanguage}[2]{{%
\expandafter\ifx\csname l@#1\endcsname\relax
\typeout{** WARNING: IEEEtran.bst: No hyphenation pattern has been}%
\typeout{** loaded for the language `#1'. Using the pattern for}%
\typeout{** the default language instead.}%
\else
\language=\csname l@#1\endcsname
\fi
#2}}
\providecommand{\BIBdecl}{\relax}
\BIBdecl

\bibitem{TSH2012}
H.~L. Trentelman, A.~A. Stoorvogel, and M.~L.~J. Hautus, \emph{Control Theory
  for Linear Systems}.\hskip 1em plus 0.5em minus 0.4em\relax Springer Science
  \& Business Media, 2012.

\bibitem{MY1979}
H.~Mayeda and T.~Yamada, ``Strong structural controllability,'' \emph{SIAM
  Journal on Control and Optimization}, vol.~17, no.~1, pp. 123--138, 1979.

\bibitem{B1981}
W.~Bachmann, ``Strenge strukturelle steuerbarkeit und beobachtbarkeit von
  mehrgr{\"o}{\ss}ensystemen/strong structural controllability and
  observability of multi-variable systems,'' \emph{at-Automatisierungstechnik},
  vol.~29, no. 1-12, pp. 318--323, 1981.

\bibitem{RSW1992}
K.~J. Reinschke, F.~Svaricek, and H.-D. Wend, ``On strong structural
  controllability of linear systems,'' in \emph{Proc. of the IEEE Conference on
  Decision and Control (CDC)}, vol.~1, 1992, pp. 203--208.

\bibitem{JFB2011}
J.~C. Jarczyk, F.~Svaricek, and B.~Alt, ``Strong structural controllability of
  linear systems revisited,'' in \emph{Proc. of the 50th IEEE Conference on
  Decision and Control (CDC) and European Control Conference (ECC)}, 2011, pp.
  1213--1218.

\bibitem{CM2013}
A.~Chapman and M.~Mesbahi, ``On strong structural controllability of networked
  systems: A constrained matching approach,'' in \emph{Proc. of the American
  Control Conference (ACC)}, 2013, pp. 6126--6131.

\bibitem{MZC2014}
N.~Monshizadeh, S.~Zhang, and M.~K. Camlibel, ``Zero forcing sets and
  controllability of dynamical systems defined on graphs,'' \emph{IEEE
  Transactions on Automatic Control}, vol.~59, no.~9, pp. 2562--2567, 2014.

\bibitem{TD2015}
M.~Trefois and J.-C. Delvenne, ``Zero forcing number, constrained matchings and
  strong structural controllability,'' \emph{Linear Algebra and its
  Applications}, vol. 484, pp. 199--218, 2015.

\bibitem{RM2009}
R.~Guimer{\`a} and M.~Sales-Pardo, ``Missing and spurious interactions and the
  reconstruction of complex networks,'' \emph{Proceedings of the National
  Academy of Sciences}, vol. 106, no.~52, pp. 22\,073--22\,078, 2009.

\bibitem{MHM2019}
S.~S. Mousavi, M.~Haeri, and M.~Mesbahi, ``Strong structural controllability
  under network perturbations,'' \emph{arXiv preprint arXiv:1904.09960}, 2019.

\bibitem{PPKAI2019}
N.~Popli, S.~Pequito, S.~Kar, A.~P. Aguiar, and M.~Ili\'c, ``Selective strong
  structural minimum-cost resilient co-design for regular descriptor linear
  systems,'' \emph{Automatica}, vol. 102, pp. 80--85, 2019.

\bibitem{JHHK2019}
J.~Jia, H.~J. van Waarde, H.~L. Trentelman, and M.~K. Camlibel, ``A unifying
  framework for strong structural controllability,'' \emph{to appear in IEEE
  Transactions on Automatic Control}, 2021.

\bibitem{MHM2017}
S.~S. Mousavi, M.~Haeri, and M.~Mesbahi, ``On the structural and strong
  structural controllability of undirected networks,'' \emph{IEEE Transactions
  on Automatic Control}, vol.~63, no.~7, pp. 2234--2241, 2017.

\bibitem{TVDF2018}
T.~Menara, V.~Katewa, D.~S. Bassett, and F.~Pasqualetti, ``The structured
  controllability radius of symmetric (brain) networks,'' in \emph{Proc. of the
  American Control Conference (ACC)}, 2018, pp. 2802--2807.

\bibitem{TDF2018}
T.~{Menara}, D.~S. {Bassett}, and F.~{Pasqualetti}, ``Structural
  controllability of symmetric networks,'' \emph{IEEE Transactions on Automatic
  Control}, vol.~64, no.~9, pp. 3740--3747, 2019.

\bibitem{LM2019}
F.~Liu and A.~S. Morse, ``A graphical characterization of structurally
  controllable linear systems with dependent parameters,'' \emph{IEEE
  Transactions on Automatic Control}, vol.~64, no.~11, pp. 4484--4495, 2019.

\bibitem{JTBC20181}
J.~Jia, H.~L. Trentelman, W.~Baar, and M.~K. Camlibel, ``A sufficient condition
  for colored strong structural controllability of networks,''
  \emph{IFAC-PapersOnLine}, vol.~51, no.~23, pp. 16--21, 2018.

\bibitem{JTBC2018}
------, ``Strong structural controllability of systems on colored graphs,''
  \emph{to appear in IEEE Transactions on Automatic Control}, 2020.

\bibitem{AIM2008}
F.~Barioli, W.~Barrett, S.~M. Fallat, H.~T. Hall, L.~Hogben, B.~Shader,
  P.~van~den Driessche, and H.~van~der Holst, ``Zero forcing parameters and
  minimum rank problems,'' \emph{Linear Algebra and its Applications}, vol.
  433, no.~2, pp. 401--411, 2010.

\bibitem{WRS2014}
A.~Weber, G.~Reissig, and F.~Svaricek, ``A linear time algorithm to verify
  strong structural controllability,'' in \emph{Proc. of the 53rd IEEE
  Conference on Decision and Control (CDC)}.\hskip 1em plus 0.5em minus
  0.4em\relax IEEE, 2014, pp. 5574--5580.

\bibitem{vWCT2017}
H.~J. van Waarde, M.~K. Camlibel, and H.~L. Trentelman, ``A distance-based
  approach to strong target control of dynamical networks,'' \emph{IEEE
  Transactions on Automatic Control}, vol.~62, no.~12, pp. 6266--6277, 2017.

\bibitem{MCT2015}
N.~Monshizadeh, M.~K. Camlibel, and H.~L.Trentelman, ``Strong targeted
  controllability of dynamical networks,'' in \emph{Proc. of the IEEE
  Conference on Decision and Control (CDC)}, 2015, pp. 4782--4787.

\bibitem{WTC2018}
H.~J. van Waarde, P.~Tesi, and M.~K. Camlibel, ``Identifiability of undirected
  dynamical networks: A graph-theoretic approach,'' \emph{IEEE Control Systems
  Letters}, vol.~2, no.~4, pp. 683--688, 2018.

\end{thebibliography}
\end{document}